\documentclass[12pt,thmsa]{article}
\usepackage{amssymb}

\usepackage{graphicx}
\usepackage{graphicx}

\begin{document}

\author{Francis OGER}
\title{\textbf{Equivalence \'{e}l\'{e}mentaire entre pavages}}
\date{}

\begin{center}
\textbf{Coverings of the plane by self-avoiding curves}

\textbf{which satisfy the local isomorphism property}

Francis OGER\bigskip \medskip
\end{center}

\noindent \textbf{Abstract.} A self-avoiding plane-filling curve cannot be
periodic, but we show that it can satisfy the local isomorphism property. We
investigate three families of coverings of the plane by finite sets of
nonoverlapping self-avoiding curves which satisfy that property in a strong
form. These curves are respectively inductive limits of: 1) $n$-folding
square curves such as the dragon curve, obtained by folding $n$ times a
strip of paper in $2$ and unfolding it with $\pi /2$ angles; 2) $n$-folding
triangular curves such as the terdragon curve, obtained by folding $n$ times
a strip of paper in $3$ and unfolding it with $\pi /3$ angles; 3)
generalizations of Peano-Gosper curves. In each family, the coverings
consist of a small number of curves (at most $6$), and in many examples only 
$1$\ curve. We do not know presently if similar examples exist in $%
%TCIMACRO{\U{211d} }%
%BeginExpansion
\mathbb{R}
%EndExpansion
^{n}$\ for $n\geq 3$.\bigskip


\bigskip

We denote by $%
%TCIMACRO{\U{2115} }%
%BeginExpansion
\mathbb{N}
%EndExpansion
^{\ast }$ the set of strictly positive integers, and $%
%TCIMACRO{\U{211d} }%
%BeginExpansion
\mathbb{R}
%EndExpansion
^{+}$\ the set of positive real numbers.

For each $n\in 
%TCIMACRO{\U{2115} }%
%BeginExpansion
\mathbb{N}
%EndExpansion
^{\ast }$, we consider $%
%TCIMACRO{\U{211d} }%
%BeginExpansion
\mathbb{R}
%EndExpansion
^{n}$\ equipped with a norm $x\rightarrow \left\| x\right\| $. For any $%
x,y\in 
%TCIMACRO{\U{211d} }%
%BeginExpansion
\mathbb{R}
%EndExpansion
^{n}$, we write $d(x,y)=\left\| y-x\right\| $. For each $x\in 
%TCIMACRO{\U{211d} }%
%BeginExpansion
\mathbb{R}
%EndExpansion
^{n}$ and each $r\in 
%TCIMACRO{\U{211d} }%
%BeginExpansion
\mathbb{R}
%EndExpansion
^{+}$, we denote by $B(x,r)$ the closed ball of center $x$ and radius $r$.

For any $E,F\subset 
%TCIMACRO{\U{211d} }%
%BeginExpansion
\mathbb{R}
%EndExpansion
^{n}$, an \emph{isomorphism} from $E$ to $F$ is a translation $\tau $ such
that $\tau (E)=F$.\ They are \emph{locally isomorphic} if, for each $x\in 
%TCIMACRO{\U{211d} }%
%BeginExpansion
\mathbb{R}
%EndExpansion
^{n}$\ (resp. $y\in 
%TCIMACRO{\U{211d} }%
%BeginExpansion
\mathbb{R}
%EndExpansion
^{n}$) and each $r\in 
%TCIMACRO{\U{211d} }%
%BeginExpansion
\mathbb{R}
%EndExpansion
^{+}$, there exists $y\in 
%TCIMACRO{\U{211d} }%
%BeginExpansion
\mathbb{R}
%EndExpansion
^{n}$\ (resp. $x\in 
%TCIMACRO{\U{211d} }%
%BeginExpansion
\mathbb{R}
%EndExpansion
^{n}$) such that $(B(x,r)\cap E,x)\cong (B(y,r)\cap F,y)$.

$E\subset 
%TCIMACRO{\U{211d} }%
%BeginExpansion
\mathbb{R}
%EndExpansion
^{n}$ satisfies the \emph{local isomorphism property} if, for each $x\in 
%TCIMACRO{\U{211d} }%
%BeginExpansion
\mathbb{R}
%EndExpansion
^{n}$ and each $r\in 
%TCIMACRO{\U{211d} }%
%BeginExpansion
\mathbb{R}
%EndExpansion
^{+}$, there exist $s\in 
%TCIMACRO{\U{211d} }%
%BeginExpansion
\mathbb{R}
%EndExpansion
^{+}$ such that each $B(y,s)$ contains some $z$ with $(B(z,r)\cap E,z)\cong
(B(x,r)\cap E,x)$. We say that $E$ satisfies the \emph{strong local
isomorphism property} if there exist $s,t\in 
%TCIMACRO{\U{211d} }%
%BeginExpansion
\mathbb{R}
%EndExpansion
^{+}$ such that, for each $x\in 
%TCIMACRO{\U{211d} }%
%BeginExpansion
\mathbb{R}
%EndExpansion
^{n}$ and each $r\in 
%TCIMACRO{\U{211d} }%
%BeginExpansion
\mathbb{R}
%EndExpansion
^{+}$, each $B(y,rs+t)$ contains some $z$ with $(B(z,r)\cap E,z)\cong
(B(x,r)\cap E,x)$. These definitions are similar to those for aperiodic
tilings.

For $n\geq 2$, we say that a closed subset $C\subset 
%TCIMACRO{\U{211d} }%
%BeginExpansion
\mathbb{R}
%EndExpansion
^{n}$\ is a \emph{self-avoiding curve} if there exists a bicontinuous map
from a closed interval $I\subset 
%TCIMACRO{\U{211d} }%
%BeginExpansion
\mathbb{R}
%EndExpansion
$ to $C$. If $I=[a,b]$ with $a,b\in 
%TCIMACRO{\U{211d} }%
%BeginExpansion
\mathbb{R}
%EndExpansion
$ (resp. $I=[a,+\infty \lbrack $ with $a\in 
%TCIMACRO{\U{211d} }%
%BeginExpansion
\mathbb{R}
%EndExpansion
$, $I=%
%TCIMACRO{\U{211d} }%
%BeginExpansion
\mathbb{R}
%EndExpansion
$), then $C$ is a \emph{bounded curve} (resp. a \emph{half curve}, a \emph{%
complete curve}).

We say that $A\subset 
%TCIMACRO{\U{211d} }%
%BeginExpansion
\mathbb{R}
%EndExpansion
^{n}$ is \emph{space-filling}\ if there exists $\alpha \in 
%TCIMACRO{\U{211d} }%
%BeginExpansion
\mathbb{R}
%EndExpansion
^{+}$\ such that $\sup_{x\in 
%TCIMACRO{\U{211d} }%
%BeginExpansion
\mathbb{R}
%EndExpansion
^{n}}\inf_{y\in A}d(x,y)\leq \alpha $. It follows from the proposition below
that a self-avoiding space-filling complete curve\ cannot be invariant
through a nontrivial translation:\bigskip

\noindent \textbf{Proposition.} For each integer $n\geq 2$, each
self-avoiding complete curve $C\subset 
%TCIMACRO{\U{211d} }%
%BeginExpansion
\mathbb{R}
%EndExpansion
^{n}$\ and each $w\in 
%TCIMACRO{\U{211d} }%
%BeginExpansion
\mathbb{R}
%EndExpansion
^{n}$\ such that $C=w+C$, there exists $B\subset C$\ bounded such that $C$
is the disjoint union of the subsets $kw+B$\ for $k\in 
%TCIMACRO{\U{2124} }%
%BeginExpansion
\mathbb{Z}
%EndExpansion
$.\bigskip

\noindent \textbf{Proof.} Consider $f:%
%TCIMACRO{\U{211d} }%
%BeginExpansion
\mathbb{R}
%EndExpansion
\rightarrow C$ bicontinuous.\ Then $\tau _{w}:C\rightarrow C:y\rightarrow
w+y $ and $\varphi :%
%TCIMACRO{\U{211d} }%
%BeginExpansion
\mathbb{R}
%EndExpansion
\rightarrow 
%TCIMACRO{\U{211d} }%
%BeginExpansion
\mathbb{R}
%EndExpansion
:s\rightarrow f^{-1}(w+f(s))$ satisfy $\varphi =f^{-1}\circ \tau _{w}\circ f$%
. The map $\varphi $ is bicontinuous because $f$, $\tau _{w}$\ and $f^{-1}$
are bicontinuous. In particular, $\varphi $ is strictly increasing or
strictly decreasing. Actually, $\varphi $ is strictly increasing since we
have $\varphi (s)\neq s$\ for each $s\in 
%TCIMACRO{\U{211d} }%
%BeginExpansion
\mathbb{R}
%EndExpansion
$.

Now consider $x\in C$\ and write $t_{k}=$ $f^{-1}(kw+x)$\ for each $k\in 
%TCIMACRO{\U{2124} }%
%BeginExpansion
\mathbb{Z}
%EndExpansion
$. The equalities $\varphi (t_{k})=t_{k+1}$ with $t_{k}<t_{k+1}$\ or $%
t_{k+1}<t_{k}$\ for $k\in 
%TCIMACRO{\U{2124} }%
%BeginExpansion
\mathbb{Z}
%EndExpansion
$ imply $\varphi ([t_{k},t_{k+1}[)=[t_{k+1},t_{k+2}[$ and $%
w+f([t_{k},t_{k+1}[)=f([t_{k+1},t_{k+2}[)$ for each $k\in 
%TCIMACRO{\U{2124} }%
%BeginExpansion
\mathbb{Z}
%EndExpansion
$. Consequently, $C$ is the union of the disjoint subsets $%
kw+f([t_{0},t_{1}[)=f([t_{k},t_{k+1}[)$\ for $k\in 
%TCIMACRO{\U{2124} }%
%BeginExpansion
\mathbb{Z}
%EndExpansion
$.~~$\blacksquare $\bigskip

On the other hand, the results of the next three sections give various
examples of self-avoiding plane-filling complete curves which satisfy the
strong local isomorphism property.

In each of them, we start with a set $\mathcal{B}$ of self-avoiding bounded
curves. We consider the set $\mathcal{C}$ of complete curves obtained as
inductive limits of curves in $\mathcal{B}$. We show that any curve in $%
\mathcal{C}$ can be extended into a generally unique plane-filling set of at
least $1$ and at most $6$ nonoverlapping such curves which satisfies the
local isomorphism property. In each case, the bounded self-avoiding curves
that we consider are described in B. Mandelbrot's book [3].

We do not know presently if similar examples exist in $%
%TCIMACRO{\U{211d} }%
%BeginExpansion
\mathbb{R}
%EndExpansion
^{n}$\ for $n\geq 3$.\bigskip

\textbf{1. Dragons and other square folding curves}\bigskip

We consider a regular tiling of the plane by squares. For each $n\in 
%TCIMACRO{\U{2115} }%
%BeginExpansion
\mathbb{N}
%EndExpansion
$, an $n$\emph{-folding curve} is obtained by folding $n$ times a strip of
paper in two, each time possibly to the left or to the right, then unfolding
it with $\pi /2$ angles so that the support of each segment of the curve is
a side of one of the squares. If all the foldings are done in the same
direction, then we obtain the dragon curve of order $n$.

For each $n$-folding curve, each side of each square is the support of at
most $1$ segment, so that we make the curve self-avoiding by rounding the
angles.

The \emph{complete folding curves} are the complete curves obtained as
inductive limits of $n$-folding curves for $n\in 
%TCIMACRO{\U{2115} }%
%BeginExpansion
\mathbb{N}
%EndExpansion
$. We say that a set of such curves is a \emph{covering of the plane} if
each side of each square is the support of exactly $1$ segment of $1$ curve.
Then the curves are disjoint.

Here, we consider nonoriented curves. The following results were proved in
[4] and [5]:\bigskip

\noindent \textbf{Theorem 1.1.} Each complete folding curve has a unique
extension into a covering of the plane by such curves which satisfies the
local isomorphism property, except in one case where there are $2$ such
extensions. Each covering actually satisfies the strong local isomorphism
property.\bigskip

\noindent \textbf{Theorem 1.2.} Each such covering consists of $1$, $2$, $3$%
, $4$ or $6$ curves.\bigskip

\noindent \textbf{Theorem 1.3.} There exist $2^{\omega }$ local isomorphism
classes of such coverings.\bigskip

We note that, for any such coverings $\mathcal{C},\mathcal{D}$ and any
curves $C\in \mathcal{C}$, $D\in \mathcal{D}$, the coverings $\mathcal{C},%
\mathcal{D}$ are locally isomorphic if and only if, for each $n\in 
%TCIMACRO{\U{2115} }%
%BeginExpansion
\mathbb{N}
%EndExpansion
$, the isomorphism classes of $n$-folding subcurves of $C$ and $D$ are the
same.\bigskip

\noindent \textbf{Theorem 1.4.} Each local isomorphism class is the union of 
$2^{\omega }$ isomorphism classes, including one with a covering by $1$
curve and one with a covering by $2$ curves.\bigskip

In [5], the local isomorphism classes which contain a covering, and actually
two nonisomorphic coverings, by $6$ curves are characterized; an example is
given by Figure 1.1 below. Some examples of classes with coverings by $3$
(resp. $4$) curves are also given. It would be intersting to determine for
each local isomorphism class the number of isomorphism classes of coverings
by $1$, $2$, $3$ or $4$ curves, especially since similar results are
obtained in the two next sections for triangular folding curves and for
Peano-Gosper curves.

\begin{center}
\includegraphics[scale=0.52]{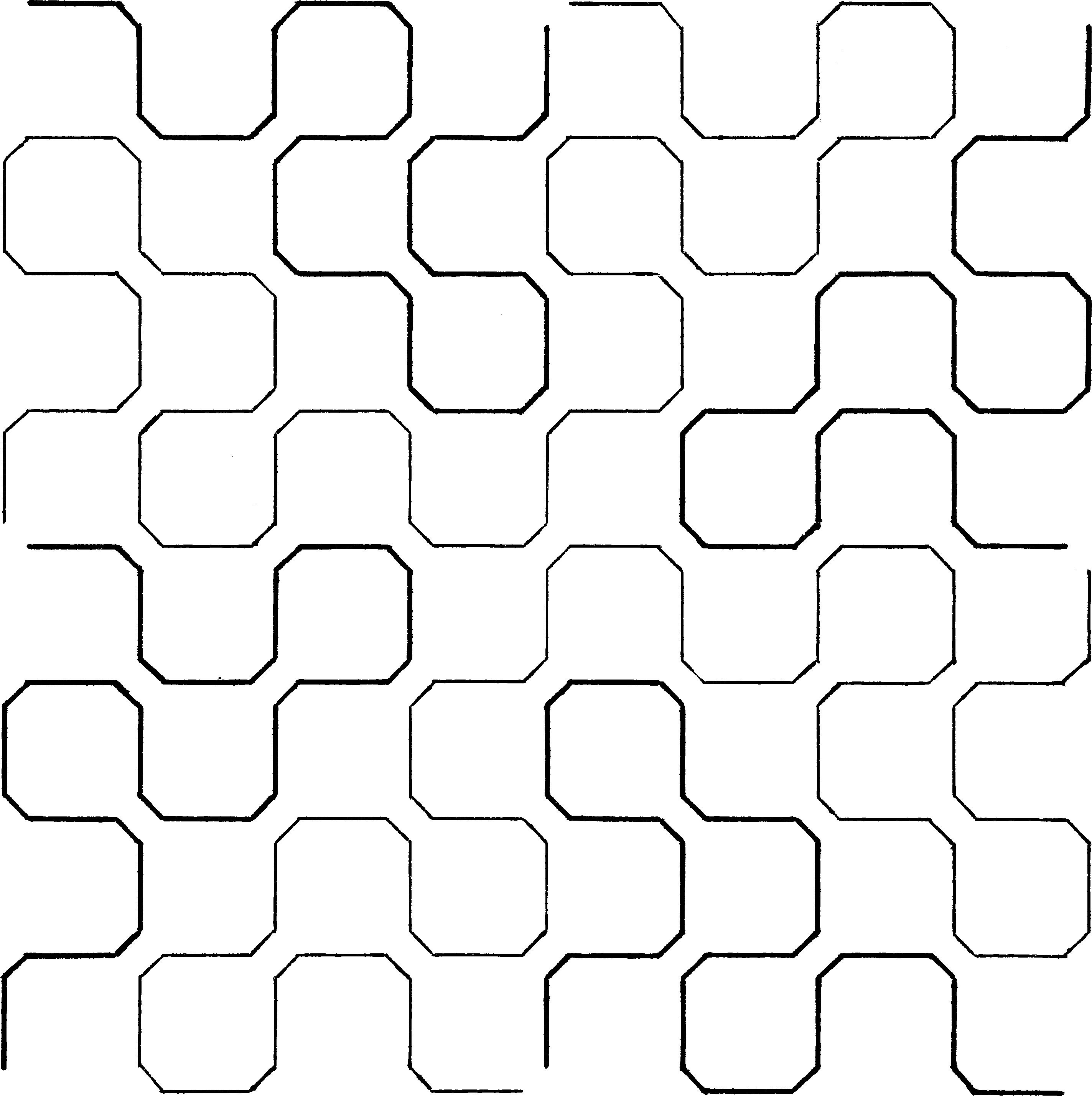}

\medskip Figure 1.1
\end{center}

\textbf{2. Terdragons and other triangular folding curves}\bigskip

For each $n\in 
%TCIMACRO{\U{2115} }%
%BeginExpansion
\mathbb{N}
%EndExpansion
$, we consider $n$-folding triangular curves such as the terdragon curve,
which are obtained by folding $n$ times a strip of paper in $3$, each time
possibly left then right or right then left, and unfolding it with $\pi /3$
angles. We prove that they are self-avoiding.

The complete folding triangular curves are obtained as inductive limits of $%
n $-folding triangular curves for $n\in 
%TCIMACRO{\U{2115} }%
%BeginExpansion
\mathbb{N}
%EndExpansion
$. We show that each nonoriented such curve can be extended into a unique
covering of the plane by disjoint such curves and that this covering
satisfies the strong local isomorphism property. We also prove that there
exist $2^{\omega }$ local isomorphism classes of such coverings and that
each of them contains $2^{\omega }$ (resp. $2^{\omega }$, $2$ or $5$, $0$)
isomorphism classes of coverings by $1$ (resp. $2$, $3$, $\geq 4$)
curves.\bigskip

\textbf{2.1. Definitions and main results}.\bigskip

Here, we consider a regular tiling $\mathcal{P}$ of the plane by equilateral
triangles. We define \emph{oriented triangular curves} which we call \emph{%
t-curves}. A \emph{bounded t-curve} (resp. \emph{half t-curve}, \emph{%
complete t-curve}) is a sequence $(A_{k})_{0\leq k\leq n}$ (resp. $%
(A_{k})_{k\in 
%TCIMACRO{\U{2115} }%
%BeginExpansion
\mathbb{N}
%EndExpansion
}$, $(A_{k})_{k\in 
%TCIMACRO{\U{2124} }%
%BeginExpansion
\mathbb{Z}
%EndExpansion
}$) of segments which are oriented sides of triangles of $\mathcal{P}$ such
that, for any $A_{k},A_{k+1}$, the terminal point of $A_{k}$ is the initial
point of $A_{k+1}$\ and $A_{k},A_{k+1}$\ form a $\mp \pi /3$ angle. We
associate to each such curve the sequence $(a_{k})_{1\leq k\leq n}$ (resp. $%
(a_{k})_{k\in 
%TCIMACRO{\U{2115} }%
%BeginExpansion
\mathbb{N}
%EndExpansion
^{\ast }}$, $(a_{k})_{k\in 
%TCIMACRO{\U{2124} }%
%BeginExpansion
\mathbb{Z}
%EndExpansion
}$) with $a_{k}=+1$ (resp. $a_{k}=-1$) for each $k$ such that we turn left
(resp. right) when we pass from $A_{k-1}$ to $A_{k}$.

We say that a t-curve $C$ is \emph{self-avoiding} if each side of a triangle
is the support of at most one segment of $C$. A set $\mathcal{C}$ of
t-curves is a \emph{covering of the plane} if each side of a triangle is the
support of exactly $1$ segment of $1$ curve of $\mathcal{C}$. We represent
the curves with slightly rounded angles, so that they do not pass through
the vertices of triangles. Then each self-avoiding t-curve passes at most
once through each point of the plane and the curves in a covering are
disjoint.

We define by induction on $n\in 
%TCIMACRO{\U{2115} }%
%BeginExpansion
\mathbb{N}
%EndExpansion
$ the sequences $T_{\lambda _{1}\cdots \lambda _{n}}$\ for $\lambda
_{1},\ldots ,\lambda _{n}\in \left\{ -1,+1\right\} $. We denote by $T$ the
empty sequence. For each $n\in 
%TCIMACRO{\U{2115} }%
%BeginExpansion
\mathbb{N}
%EndExpansion
$ and any $\lambda _{1},\ldots ,\lambda _{n+1}$, we write $T_{\lambda
_{1}\cdots \lambda _{n+1}}=(T_{\lambda _{1}\cdots \lambda _{n}},\lambda
_{n+1},T_{\lambda _{1}\cdots \lambda _{n}},-\lambda _{n+1},T_{\lambda
_{1}\cdots \lambda _{n}})$.

For each $n\in 
%TCIMACRO{\U{2115} }%
%BeginExpansion
\mathbb{N}
%EndExpansion
$,\ an $n$\emph{-folding t-curve} associated to $T_{\lambda _{1}\cdots
\lambda _{n}}$\ is realized as follows: We successively fold $n$ times in
three a strip of paper. For each $k\in \left\{ 1,\ldots ,n\right\} $, the $k$%
-th folding is done left then right if $\lambda _{n+1-k}=+1$ and right then
left if $\lambda _{n+1-k}=-1$. Then we unfold the strip, keeping a $\pi /3$\
angle for each fold. The \emph{terdragon} curve of order $n$, which was
first considered in [1], is the curve $T_{\lambda _{1}\cdots \lambda _{n}}$\
with $\lambda _{1}=\cdots =\lambda _{n}=+1$.

The $\infty $\emph{-folding t-curves} (resp. \emph{complete folding t-curves}%
) are the half t-curves (resp. complete t-curves) obtained as inductive
limits of $n$-folding t-curves $C_{n}$ for $n\in 
%TCIMACRO{\U{2115} }%
%BeginExpansion
\mathbb{N}
%EndExpansion
$. We are going to see that, for each complete folding t-curve $C$, there
exists a unique sequence $\Lambda =(\lambda _{n})_{n\in 
%TCIMACRO{\U{2115} }%
%BeginExpansion
\mathbb{N}
%EndExpansion
^{\ast }}\in \left\{ -1,+1\right\} ^{%
%TCIMACRO{\U{2115} }%
%BeginExpansion
\mathbb{N}
%EndExpansion
^{\ast }}$, not depending on the orientation of $C$, such that each bounded
subcurve of $C$ is countained in a t-curve associated to some $T_{\lambda
_{1}\cdots \lambda _{n}}$.

In the second part of the present section, we show that folding t-curves are
self-avoiding.

We denote by (P) the following property of a set $E$ of oriented sides of
triangles of $\mathcal{P}$: If $A,B\in E$ are sides of the same triangle,
then they define the same direction of rotation around its center.

There are $2$ opposite sets $E_{1},E_{2}$ which satisfy (P) and such that
each nonoriented side of a triangle of $\mathcal{P}$ is the support of $1$
element of $E_{1}$\ and $1$ element of $E_{2}$. The segments of any t-curve
are all contained in $E_{1}$, or all contained in $E_{2}$, and therefore
satisfy (P).

We say that a covering $\mathcal{C}$ of the plane by oriented complete
t-curves satisfies (P) if the set of segments of curves of $\mathcal{C}$
satisfies (P), or equivalently if it is equal to $E_{1}$ or $E_{2}$. Each
covering by nonoriented complete t-curves induces $2$ covering by oriented
complete t-curves which satisfy (P) and have opposite orientations.\bigskip

\noindent \textbf{Theorem 2.1.} Each nonoriented (resp. oriented) complete
folding t-curve can be extended in a unique way into a covering of the plane
by such curves (resp. a covering of the plane by such curves which satisfies
(P)). This covering satisfies the strong local isomorphism property and
consists of curves associated to the same sequence $\Lambda $.\bigskip

\noindent \textbf{Theorem 2.2.} Two coverings of the plane by nonoriented
complete folding t-curves are locally isomorphic if and only if they are
associated to the same sequence $\Lambda $.\bigskip

\noindent \textbf{Theorem 2.3.} For each covering of the plane by oriented
complete folding t-curves, the local isomorphism property implies (P) and
(P) implies the strong local isomorphism property. Any coverings with these
properties are locally isomorphic if and only if they are associated to the
same sequence $\Lambda $ and have the same orientation.\bigskip

\noindent \textbf{Theorem 2.4.} For each $\Lambda =(\lambda _{n})_{n\in 
%TCIMACRO{\U{2115} }%
%BeginExpansion
\mathbb{N}
%EndExpansion
^{\ast }}\in \left\{ -1,+1\right\} ^{%
%TCIMACRO{\U{2115} }%
%BeginExpansion
\mathbb{N}
%EndExpansion
^{\ast }}$, the class of coverings of the plane by nonoriented complete
folding t-curves associated to $\Lambda $ is the union of:

\noindent 1) $2^{\omega }$ isomorphism classes of coverings by $1$ curve;

\noindent 2) $2^{\omega }$ isomorphism classes of coverings by $2$ curves;

\noindent 3) $2$ isomorphism classes of coverings by $3$ curves having $1$
vertex in common, where each of the $3$ curves is the union of $2$ $\infty $%
-folding t-curves starting from that vertex;

\noindent 4)$\ 3$ isomorphism classes of coverings by $3$ curves with $1$
curve separating the $2$ others, if there exists $k\in 
%TCIMACRO{\U{2115} }%
%BeginExpansion
\mathbb{N}
%EndExpansion
^{\ast }$ such that $\lambda _{n+1}=-\lambda _{n}$ for each $n\geq k$
(Figure 2.1 gives an example with $k=1$ and $\lambda _{1}=+1$); no such
class otherwise.

\noindent In 3) and in 4), all the coverings are equivalent up to isometry.

\begin{center}
\includegraphics[scale=0.60]{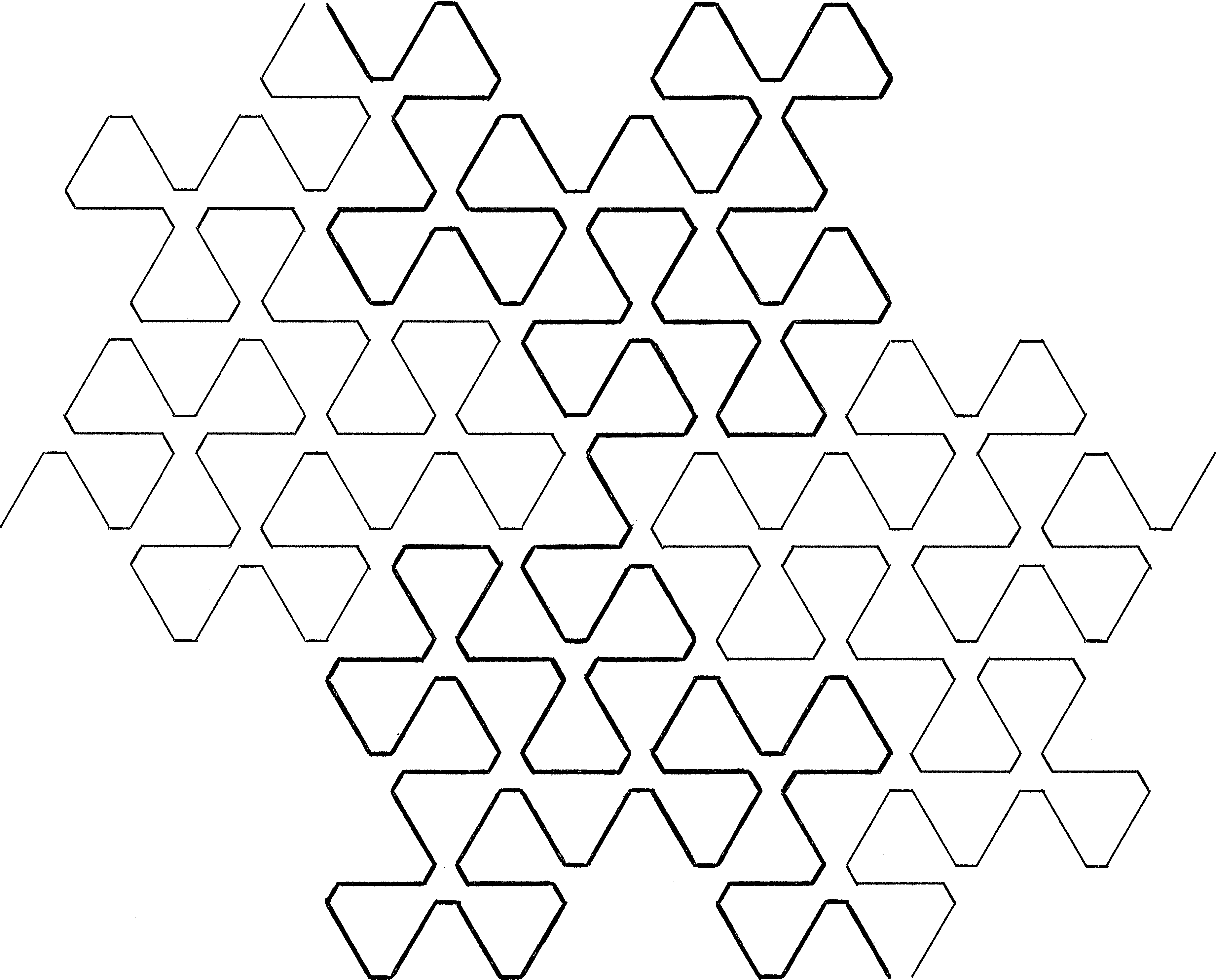}\medskip

Figure 2.1
\end{center}

\medskip

\textbf{2.2. Detailed results and proofs.}\bigskip

Unless otherwise specified, all the curves that we consider are oriented.
The property (P) and the sets $E_{1},E_{2}$ are defined as above.

For any t-curves $C,D$, if the terminal point of $C$ is the initial point of 
$D$ and if the terminal segment of $C$ and the initial segment of $D$ form a 
$\mp \pi /3$\ angle, then we denote by $CD$\ the t-curve obtained\ by
connecting them.

For each $n\in 
%TCIMACRO{\U{2115} }%
%BeginExpansion
\mathbb{N}
%EndExpansion
^{\ast }$, any $\lambda _{1},\ldots ,\lambda _{n}\in \left\{ -1,+1\right\} $
and each $n$-folding t-curve $C$ associated to $T_{\lambda _{1}\cdots
\lambda _{n}}$, we consider the $3$ $(n-1)$-folding t-curves $C^{\mathrm{I}%
},C^{\mathrm{M}},C^{\mathrm{S}}$ associated to $T_{\lambda _{1}\cdots
\lambda _{n-1}}$ such that $C=C^{\mathrm{I}}C^{\mathrm{M}}C^{\mathrm{S}}$.
Here, $\mathrm{I},\mathrm{M},\mathrm{S}$ can be viewed as abbreviations for
``inferior'', ``middle'', ``superior''.

For $k\in 
%TCIMACRO{\U{2115} }%
%BeginExpansion
\mathbb{N}
%EndExpansion
$ and $S=(\alpha _{1},\ldots ,\alpha _{k})\in \left\{ -1,+1\right\} ^{k}$,
we write $\overline{S}=(-\alpha _{k},\ldots ,-\alpha _{1})$. The reverse of
a curve associated to $S$ is associated to $\overline{S}$. For $S=(\alpha
_{n})_{n\in 
%TCIMACRO{\U{2115} }%
%BeginExpansion
\mathbb{N}
%EndExpansion
^{\ast }}\in \left\{ -1,+1\right\} ^{%
%TCIMACRO{\U{2115} }%
%BeginExpansion
\mathbb{N}
%EndExpansion
^{\ast }}$, we write $\overline{S}=(-\alpha _{-n})_{n\in -%
%TCIMACRO{\U{2115} }%
%BeginExpansion
\mathbb{N}
%EndExpansion
^{\ast }}$.

For $n\in 
%TCIMACRO{\U{2115} }%
%BeginExpansion
\mathbb{N}
%EndExpansion
$ and $\lambda _{1},\ldots ,\lambda _{n}\in \left\{ -1,+1\right\} $, we have 
$T_{\lambda _{1}\cdots \lambda _{n}}=\overline{T_{\lambda _{1}\cdots \lambda
_{n}}}$. Consequently,\ the reverse of a curve associated to $T_{\lambda
_{1}\cdots \lambda _{n}}$ is also associated to $T_{\lambda _{1}\cdots
\lambda _{n}}$.

For\ each $\Lambda =(\lambda _{k})_{k\in 
%TCIMACRO{\U{2115} }%
%BeginExpansion
\mathbb{N}
%EndExpansion
^{\ast }}\in \left\{ -1,+1\right\} ^{%
%TCIMACRO{\U{2115} }%
%BeginExpansion
\mathbb{N}
%EndExpansion
^{\ast }}$, we denote by $T_{\Lambda }$ the inductive limit of the sequences 
$T_{\lambda _{1}\cdots \lambda _{n}}$ with $T_{\lambda _{1}\cdots \lambda
_{n}}$ initial segment of $T_{\lambda _{1}\cdots \lambda _{n+1}}$ for each $%
n\in 
%TCIMACRO{\U{2115} }%
%BeginExpansion
\mathbb{N}
%EndExpansion
$.\bigskip

\noindent \textbf{Proposition 2.5.} Let $C$ be a complete t-curve associated
to a sequence $S=(s_{h})_{h\in 
%TCIMACRO{\U{2124} }%
%BeginExpansion
\mathbb{Z}
%EndExpansion
}\in \left\{ -1,+1\right\} ^{%
%TCIMACRO{\U{2124} }%
%BeginExpansion
\mathbb{Z}
%EndExpansion
}$. Then the properties 1), 2), 3) below are equivalent:

\noindent 1) $C$ is a complete folding t-curve.

\noindent 2) For each $k\in 
%TCIMACRO{\U{2115} }%
%BeginExpansion
\mathbb{N}
%EndExpansion
$, there exists $h\in 
%TCIMACRO{\U{2124} }%
%BeginExpansion
\mathbb{Z}
%EndExpansion
$ such that $s_{h+3^{k}+3^{k+1}i}=-s_{h+3^{k}2+3^{k+1}j}$\ for any $i,j\in 
%TCIMACRO{\U{2124} }%
%BeginExpansion
\mathbb{Z}
%EndExpansion
$.

\noindent 3) There exists a unique sequence $\Lambda =(\lambda _{n})_{n\in 
%TCIMACRO{\U{2115} }%
%BeginExpansion
\mathbb{N}
%EndExpansion
^{\ast }}\in \left\{ -1,+1\right\} ^{%
%TCIMACRO{\U{2115} }%
%BeginExpansion
\mathbb{N}
%EndExpansion
^{\ast }}$ such that exactly one of the two following properties is true:

\noindent a) $S$ is equivalent to $(\overline{T_{\Lambda }},+1,T_{\Lambda })$
or $(\overline{T_{\Lambda }},-1,T_{\Lambda })$ modulo a translation of $%
%TCIMACRO{\U{2124} }%
%BeginExpansion
\mathbb{Z}
%EndExpansion
$;

\noindent b) $C=\cup _{n\in 
%TCIMACRO{\U{2115} }%
%BeginExpansion
\mathbb{N}
%EndExpansion
}C_{n}$\ for a sequence $(C_{n})_{n\in 
%TCIMACRO{\U{2115} }%
%BeginExpansion
\mathbb{N}
%EndExpansion
}$ such that, for each $n\in 
%TCIMACRO{\U{2115} }%
%BeginExpansion
\mathbb{N}
%EndExpansion
$, $C_{n}$\ is an $n$-folding t-curve\ associated to $T_{\lambda _{1}\cdots
\lambda _{n}}$\ and $C_{n}\in \left\{
C_{n+1}^{I},C_{n+1}^{M},C_{n+1}^{S}\right\} $.\bigskip

\noindent \textbf{Proof.} 1) or 3) implies 2) because we have $%
s_{3^{k}+3^{k+1}i}=-s_{3^{k}2+3^{k+1}j}$\ for each $n\in 
%TCIMACRO{\U{2115} }%
%BeginExpansion
\mathbb{N}
%EndExpansion
$, each sequence $(s_{h})_{1\leq h\leq 3^{n}-1}$\ associated to an $n$%
-folding t-curve, each $k\in \left\{ 0,\ldots ,n-1\right\} $\ and any
integers $0\leq i,j\leq 3^{n-k-1}-1$. If the case a) of 3) is realized,
then, in 2), we can take the same $h$\ for each $k\in 
%TCIMACRO{\U{2115} }%
%BeginExpansion
\mathbb{N}
%EndExpansion
$.

Now we show that 2) implies 1) and 3). For each $g\in 
%TCIMACRO{\U{2124} }%
%BeginExpansion
\mathbb{Z}
%EndExpansion
$, we consider the point $z_{g}$\ of the plane which is associated to $s_{g}$%
\ in the correspondance between $S$ and $C$.\ For $g\in 
%TCIMACRO{\U{2124} }%
%BeginExpansion
\mathbb{Z}
%EndExpansion
$ and $n\in 
%TCIMACRO{\U{2115} }%
%BeginExpansion
\mathbb{N}
%EndExpansion
$, we denote by $\alpha _{n}(g)$\ the largest integer $h\leq g$\ such that $%
s_{h+3^{k}+3^{k+1}i}=-s_{h+3^{k}2+3^{k+1}j}$\ for each $k\in \left\{
0,\ldots ,n-1\right\} $\ and any $i,j\in 
%TCIMACRO{\U{2124} }%
%BeginExpansion
\mathbb{Z}
%EndExpansion
$. We have $g<\beta _{n}(g)$ for $\beta _{n}(g)=\alpha _{n}(g)+3^{n}$. The
part $C_{n}(g)$\ of $C$ between $z_{\alpha _{n}(g)}$\ and $z_{\beta _{n}(g)}$%
\ is an $n$-folding t-curve. We have $C_{n}(g)\in
\{C_{n+1}(g)^{I},C_{n+1}(g)^{M},C_{n+1}(g)^{S}\}$.

There exists a unique sequence $\Lambda =(\lambda _{n})_{n\in 
%TCIMACRO{\U{2115} }%
%BeginExpansion
\mathbb{N}
%EndExpansion
^{\ast }}\in \left\{ -1,+1\right\} ^{%
%TCIMACRO{\U{2115} }%
%BeginExpansion
\mathbb{N}
%EndExpansion
^{\ast }}$ such that,\ for each $n\in 
%TCIMACRO{\U{2115} }%
%BeginExpansion
\mathbb{N}
%EndExpansion
$ and each $g\in 
%TCIMACRO{\U{2124} }%
%BeginExpansion
\mathbb{Z}
%EndExpansion
$, $C_{n}(g)$\ is associated to $T_{\lambda _{1}\cdots \lambda _{n}}$. For
any $g,h\in 
%TCIMACRO{\U{2124} }%
%BeginExpansion
\mathbb{Z}
%EndExpansion
$ such that $g<h$, we have $C_{n}(g)=C_{n}(h)$ for $n$ large enough, or $%
\beta _{n}(g)=\alpha _{n}(h)$ for $n$ large enough.

If there exists $g\in 
%TCIMACRO{\U{2124} }%
%BeginExpansion
\mathbb{Z}
%EndExpansion
$ such that $\cup _{n\in 
%TCIMACRO{\U{2115} }%
%BeginExpansion
\mathbb{N}
%EndExpansion
}C_{n}(g)=C$, then the case b) of 3) is realized and $C$ is clearly a
complete folding t-curve.

Otherwise, there exist $m\in 
%TCIMACRO{\U{2115} }%
%BeginExpansion
\mathbb{N}
%EndExpansion
^{\ast }$ and $g,h\in 
%TCIMACRO{\U{2124} }%
%BeginExpansion
\mathbb{Z}
%EndExpansion
$ with $g<h$ such that $\beta _{n}(g)=\alpha _{n}(h)$ for each $n\geq m$.
Then the case a) of 3) is realized. Moreover, for each $n\geq m$, $%
C_{n}(g)^{M}C_{n}(g)^{S}C_{n}(h)^{I}$\ or $%
C_{n}(g)^{S}C_{n}(h)^{I}C_{n}(h)^{M}$\ is an $n$-folding t-curve $D_{n}$
contained in $C$. We have $D_{n}\subset D_{n+1}$ for each $n\geq m$ and $%
C=\cup _{n\geq m}D_{n}$, whence 1).~~$\blacksquare $\bigskip

Similarly to the case of folding curves in [4] and [5], and Peano-Gosper
curves in Section 3 below, we define a \emph{derivation} $\Delta $\ on
folding t-curves.

For each $n\in 
%TCIMACRO{\U{2115} }%
%BeginExpansion
\mathbb{N}
%EndExpansion
^{\ast }$, we divide each $n$-folding t-curve $C$ into sequences of $3$\
consecutive segments; $\Delta (C)$ is obtained by replacing each such
sequence with a unique segment whose initial and terminal points are the
initial point of the first segment and the terminal point of the third
segment; if $C$ is associated to a sequence $(\alpha _{i})_{1\leq i\leq
3^{n}-1}\in \left\{ -1,+1\right\} ^{3^{n}-1}$, then $\Delta (C)$ is
associated to $(\alpha _{3i})_{1\leq i\leq 3^{n-1}-1}$.

The definition of $\Delta $ naturally extends to $\infty $-folding t-curves.
Now we extend it to complete folding t-curves.

Consider any such curve $C=$ $(A_{k})_{k\in 
%TCIMACRO{\U{2124} }%
%BeginExpansion
\mathbb{Z}
%EndExpansion
}$, and the associated sequence $(s_{k})_{k\in 
%TCIMACRO{\U{2124} }%
%BeginExpansion
\mathbb{Z}
%EndExpansion
}\in \left\{ -1,+1\right\} ^{%
%TCIMACRO{\U{2124} }%
%BeginExpansion
\mathbb{Z}
%EndExpansion
}$. Then, by Proposition 2.5, there exist $\varepsilon \in \left\{
-1,+1\right\} $ and $h\in 
%TCIMACRO{\U{2124} }%
%BeginExpansion
\mathbb{Z}
%EndExpansion
$, whose remainder modulo $3$ is completely determined by $C$, such that $%
s_{h+3k+1}=\varepsilon $ and $s_{h+3k+2}=-\varepsilon $ for each $k\in 
%TCIMACRO{\U{2124} }%
%BeginExpansion
\mathbb{Z}
%EndExpansion
$.\ The curve $\Delta (C)$ is obtained by replacing each $%
A_{h+3k}A_{h+3k+1}A_{h+3k+2}$ with a unique segment; it is associated to the
sequence $(s_{h+3k})_{k\in 
%TCIMACRO{\U{2124} }%
%BeginExpansion
\mathbb{Z}
%EndExpansion
}$.

If an $\infty $-folding or a complete folding t-curve $C$ is the inductive
limit for $n\in 
%TCIMACRO{\U{2115} }%
%BeginExpansion
\mathbb{N}
%EndExpansion
$ of some $n$-folding t-curves $C_{n}\subset C$, then $\Delta (C)$ is the
inductive limit of the curves $\Delta (C_{n})$ for $n\in 
%TCIMACRO{\U{2115} }%
%BeginExpansion
\mathbb{N}
%EndExpansion
^{\ast }$.

For each $k\in 
%TCIMACRO{\U{2115} }%
%BeginExpansion
\mathbb{N}
%EndExpansion
^{\ast }$ and each t-curve $C$ which is $n$-folding for some $n\geq k$, $%
\infty $-folding or complete folding, $\Delta ^{k}(C)$ is obtained in the
same way by replacing sequences of $3^{k}$\ consecutive segments with $1$
segment.

Now we begin to show the self-avoiding and plane-filling properties of
folding t-curves. We denote by $U$ the set of vertices of triangles of $%
\mathcal{P}$ and we fix $i\in \left\{ 1,2\right\} $.

We prove by induction on $n$ that, for each $n\in 
%TCIMACRO{\U{2115} }%
%BeginExpansion
\mathbb{N}
%EndExpansion
$, each $\Lambda \in \left\{ -1,+1\right\} ^{n}$ and each $x\in U$, there
exists a unique covering $\mathcal{C}(\Lambda ,x)$ of the plane by $n$%
-folding t-curves associated to $T_{\Lambda }$ and with segments in $E_{i}$
such that $x$ is an endpoint of some of the curves.

If $n=0$, then $\Lambda $\ is the empty sequence and we write\ $\mathcal{C}%
(\Lambda ,x)=E_{i}$. It remains to be proved that, if the property is true
for an integer $n$, then it is true for $n+1$.

If $X$ is the set of vertices of equilateral triangles which form a regular
tiling of the plane, then, for each $x\in X$, there exists a unique
partition $X=G(X,x)\cup H(X,x)$\ such that:

\noindent 1)\ $H(X,x)$ is the set of vertices of hexagons which form a
regular tiling of the plane;

\noindent 2) $x\in G(X,x)$\ and the elements of $G(X,x)$\ are the centers of
the hexagons.

\noindent Moreover, $G(X,x)$\ is also the set of vertices of equilateral
triangles which form a regular tiling of the plane.

For each $x\in U$, we\ write $V_{0}(x)=U$\ and, for each $k\in 
%TCIMACRO{\U{2115} }%
%BeginExpansion
\mathbb{N}
%EndExpansion
$, $V_{k+1}(x)=G(V_{k}(x),x)$\ and $W_{k+1}(x)=H(V_{k}(x),x)$.

For each $\Lambda \in \left\{ -1,+1\right\} ^{n}$ and each $\lambda \in
\left\{ -1,+1\right\} $, supposing that $\mathcal{C}(\Lambda ,x)$ is already
defined, we define $\mathcal{C}((\Lambda ,\lambda ),x)$ as follows: If $%
\lambda =+1$ (resp. $-1$), then, for each $y\in V_{n+1}(x)$\ and each curve $%
A\in \mathcal{C}(\Lambda ,x)$ starting from $y$, we put in $\mathcal{C}%
((\Lambda ,\lambda ),x)$ the curve $ABC$, where:

\noindent 1) $B$ is the curve of $\mathcal{C}(\Lambda ,x)$ starting from the
endpoint of $A$ and such that its first segment is just at the left (resp.
right) of the last segment of $A$;

\noindent 2) $C$ is the curve of $\mathcal{C}(\Lambda ,x)$ starting from the
endpoint of $B$ and such that its first segment is just at the right (resp.
left) of the last segment of $B$.

The following properties are true for each $n\in 
%TCIMACRO{\U{2115} }%
%BeginExpansion
\mathbb{N}
%EndExpansion
$, each $\Lambda \in \left\{ -1,+1\right\} ^{n}$ and each $x\in U$. For $%
n\geq 2$, they are proved by using $\Delta ^{n-1}$.

Each curve of $\mathcal{C}(\Lambda ,x)$ connects a pair of points of $%
V_{n}(x)$ with minimal distance. Each such pair $(y,z)$\ is connected by a
unique curve of $\mathcal{C}(\Lambda ,x)$. If $n\geq 1$, then this curve
contains the $2$\ points of $W_{n}(x)$ which are between $y$ and $z$.

Using the derivations, we see that, if $C$ is a curve associated to some $%
T_{\Lambda }$, then $C$ or its reverse belongs to $\mathcal{C}(\Lambda ,x)=%
\mathcal{C}(\Lambda ,y)$, where $x$ and $y$ are the endpoints of $C$. It
follows that $C$ is self-avoiding.

For each $n\in 
%TCIMACRO{\U{2115} }%
%BeginExpansion
\mathbb{N}
%EndExpansion
^{\ast }$, we call an $n$\emph{-triangle} any equilateral triangle such that
each side consists of $n$ sides of triangles of $\mathcal{P}$.\ We say that
a set of curves $\mathcal{C}$ \emph{covers} an $n$-triangle $\mathcal{T}$
if, in each $1$-triangle of $\mathcal{P}$ contained in $\mathcal{T}$, each
side is the support of a segment of a curve of $\mathcal{C}$.\bigskip

\noindent \textbf{Proposition 2.6.} For each integer $n\geq 2$, each $(2n)$%
-folding t-curve covers a $(3^{n-1}+3)/2$-triangle.

\begin{center}
\includegraphics[scale=0.60]{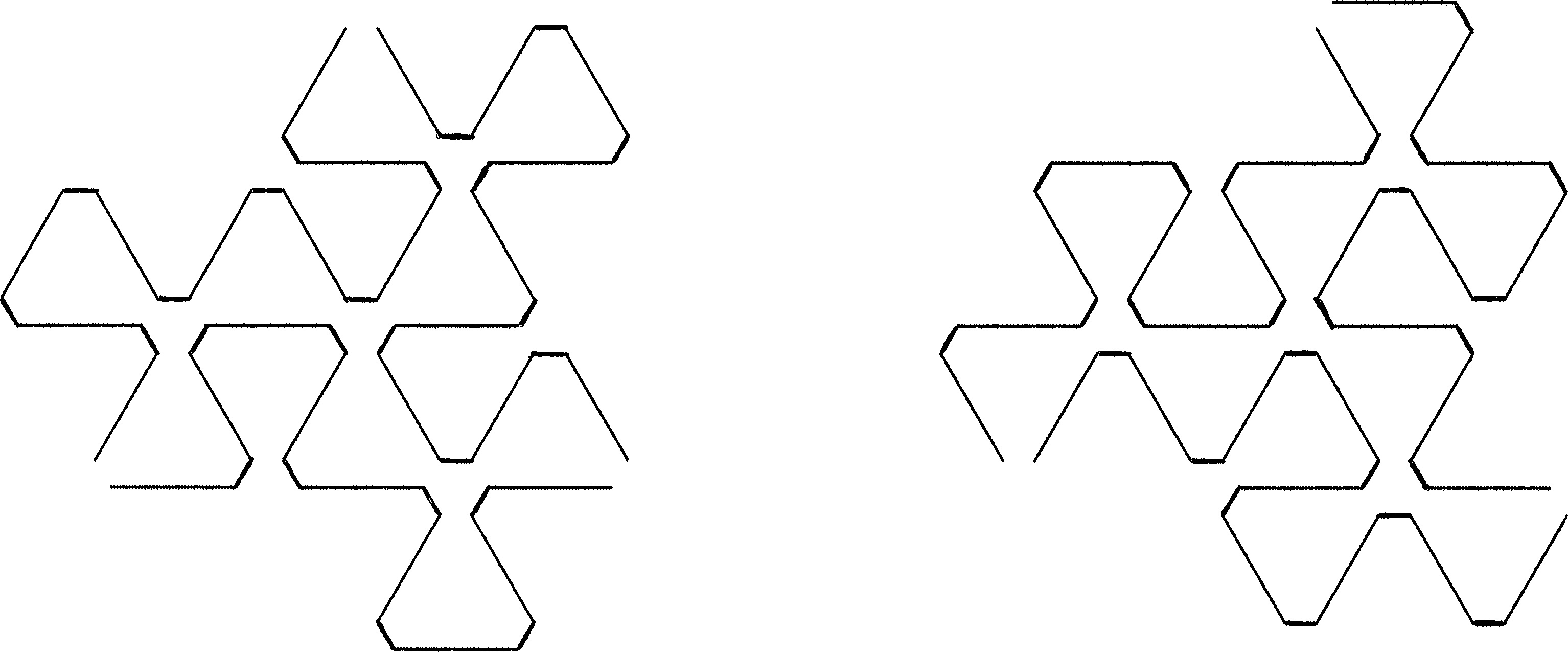}\medskip

Figure 2.2
\end{center}

\noindent \textbf{Proof.} For each $n\geq 2$, we denote by $k_{n}$\ the
largest integer $k$\ such that each $(2n)$-folding t-curve covers a $k$%
-triangle. We simultaneously prove that $k_{2}\geq 3$\ and $k_{n}\geq
3k_{n-1}-3$ for each $n\geq 3$.

For each $(2n)$-folding t-curve $C$, we consider the tiling $\mathcal{Q}$\
of the plane by equilateral triangles which is associated to $\Delta ^{2}(C)$%
. Each triangle of $\mathcal{Q}$ is the union of $9$ triangles of $\mathcal{P%
}$.

For any segments $A_{1},A_{2},A_{3}$ of $\Delta ^{2}(C)$, if their supports
are the sides of a $1$-triangle $\mathcal{T}$ of $\mathcal{Q}$, then their
orientations\ define the same direction of rotation around the center of $%
\mathcal{T}$ since $\Delta ^{2}(C)$ satisfies (P). Figure 2.2 above shows
two possible configurations for $\Delta ^{-2}(A_{1})$, $\Delta ^{-2}(A_{2})$%
, $\Delta ^{-2}(A_{3})$. Any other configuration is equivalent to one of
them up to isometry.

We see from Figure 2.2 that $\mathcal{E}=\left\{ \Delta ^{-2}(A_{1}),\Delta
^{-2}(A_{2}),\Delta ^{-2}(A_{3})\right\} $ necessarily covers a $3$-triangle
of $\mathcal{P}$. If $n=2$, then $C$ covers a $3$-triangle of $\mathcal{P}$
since the $2$-folding t-curve $\Delta ^{2}(C)$ has such segments $%
A_{1},A_{2},A_{3}$.

We also see from Figure 2.2 that each side of a $1$-triangle of $\mathcal{P}$
contained in $\mathcal{T}$ is the support of a segment of $\mathcal{E}$,
except possibly if one of its endpoints is a vertex of $\mathcal{T}$.

In order to prove that $k_{n}\geq 3k_{n-1}-3$ for each $n\geq 3$, it
suffices to show that, for each integer $k\geq 2$, if $\Delta ^{2}(C)$\
covers a $k$-triangle $\mathcal{W}\ $of $\mathcal{Q}$, then $C$ covers the $%
(3k-3)$-triangle $\mathcal{V}$ of $\mathcal{P}$ contained in the interior of 
$\mathcal{W}$.

This property is a consequence of the two following facts: First, by the
argument above, for each $1$-triangle $\mathcal{T}$ of $\mathcal{Q}$, if $%
\Delta ^{2}(C)$ covers $\mathcal{T}$, then, in each $1$-triangle of $%
\mathcal{P}$ contained in $\mathcal{T}$, each side is the support of a
segment of $C$ if neither of its endpoints is a vertex of $\mathcal{T}$.
Second, for each vertex $z$ of a $1$-triangle of $\mathcal{Q}$, if $z$
belongs to $\mathcal{V}$, then $z$ is an endpoint of $6$ segments of $C$
since it is an endpoint of $6$ segments of $\Delta ^{2}(C)$; as $C$ is
self-avoiding, it follows that\ each side with endpoint $z$ of a $1$%
-triangle of $\mathcal{P}$ is the support of a segment of $C$.~~$%
\blacksquare $\bigskip

\noindent \textbf{Corollary 2.7.} For each covering of the plane by complete
folding t-curves, the local isomorphism property implies (P).\bigskip

\noindent \textbf{Proof.} Let $\mathcal{C}$ be such a covering. By
Proposition 2.6, each $C\in \mathcal{C}$\ covers arbitrarily large
triangles. The covering $\mathcal{C}$ satisfies (P) on any such triangle,
and therefore on the whole plane by the local isomorphism property.~~$%
\blacksquare $\bigskip

For each nontrivial t-curve $C$, we say that $x\in U$ belongs to $C$ and we
write $x\in C$ if $x$ is an endpoint of at least one segment of $C$. For any 
$x,y\in U$ with $d(x,y)=1$, we say that the nonoriented segment $\left[ x,y%
\right] $ belongs to $C$ and we write $\left[ x,y\right] \in C$ if $\left[
x,y\right] $ is the support of a segment of $C$.

For each folding t-curve $C$, we denote by $\mathrm{F}_{\mathrm{L}}(C)$
(resp. $\mathrm{F}_{\mathrm{R}}(C)$) the union of the sides $\left[ x,y%
\right] $ of triangles of $\mathcal{P}$ such that $x,y\in C$, such that $C$
contains exactly $1$ of the $2$ points which form equilateral triangles with 
$x$ and $y$, and such that the second point is at the left (resp. right) of $%
C$.

For each integer $n\geq 1$ and each $n$-folding t-curve $C$\ with initial
point $u$ and terminal point $v$, there exists a unique sequence $%
(y_{i})_{1\leq i\leq h}\subset U$\ with $y_{1}=u$, $y_{h}=v$, $%
d(y_{i-1},y_{i})=1$\ and $[y_{i-1},y_{i}]\subset \mathrm{F}_{\mathrm{L}}(C)$
for $2\leq i\leq h$. There is also a unique sequence $(z_{j})_{1\leq j\leq
k}\subset U$\ with $z_{1}=u$, $z_{k}=v$, $d(z_{j-1},z_{j})=1$\ and $%
[z_{j-1},z_{j}]\subset \mathrm{F}_{\mathrm{R}}(C)$ for $2\leq j\leq k$.

We denote by $\mathrm{F}_{\mathrm{L}}^{0}(C)$ (resp. $\mathrm{F}_{\mathrm{R}%
}^{0}(C)$) the union of the segments $[y_{i-1},y_{i}]$ (resp. $%
[z_{j-1},z_{j}]$). We have $\mathrm{F}_{\mathrm{L}}^{0}(C)\cap \mathrm{F}_{%
\mathrm{R}}^{0}(C)=\left\{ u,v\right\} $ and no point of $C$ is outside $%
\mathrm{F}_{\mathrm{L}}^{0}(C)\cup \mathrm{F}_{\mathrm{R}}^{0}(C)$.

It follows that $C$ starts at $y_{1}$, then successively passes through $%
y_{2},\ldots ,y_{h-1}$ in that order, possibly once or twice for each of
them, then ends at $y_{h}$. The same property is true for $z_{1},\ldots
,z_{k}$.

Now we show that $C$ contains the segments $\left[ w_{1},w_{2}\right] $\
with $w_{1},w_{2}\in U$\ and $d(w_{1},w_{2})=1$\ which are inside $\mathrm{F}%
_{\mathrm{L}}^{0}(C)\cup \mathrm{F}_{\mathrm{R}}^{0}(C)$. It follows $%
\mathrm{F}_{\mathrm{L}}(C)=\mathrm{F}_{\mathrm{L}}^{0}(C)$\ and $\mathrm{F}_{%
\mathrm{R}}(C)=\mathrm{F}_{\mathrm{R}}^{0}(C)$.

For the proof, we use the existence of a covering $\mathcal{C}(\Lambda ,u)$
of the plane which contains $C$.\ If a segment $\left[ w_{1},w_{2}\right] $\
as above does not belong to $C$, then it belongs to a curve $D\in \mathcal{C}%
(\Lambda ,u)$ with $D\neq C$.

The segment $\left[ w_{1},w_{2}\right] $ is between $\mathrm{F}_{\mathrm{L}%
}^{0}(C)$ and $C$, or between $C$ and $\mathrm{F}_{\mathrm{R}}^{0}(C)$. We
can suppose that the first case is true, since the second one can be treated
in the same way.

Then there exists $i$ such that $\left[ w_{1},w_{2}\right] $ is inside the
loop formed by $C$ between two occurrences of $x_{i}$, or between an
occurrence of $x_{i}$ and an occurrence of $x_{i+1}$. We observe that, in
the first case, no curve of $\mathcal{C}(\Lambda ,u)$ crosses the limit of
the loop, and in the second case at most one curve of $\mathcal{C}(\Lambda
,u)$ crosses that limit and only once, since the part of the curve which
crosses the limit must contain $\left[ x_{i},x_{i+1}\right] $.

In both cases, it follows that at least one endpoint of $D$ is inside the
loop. Then another curve of $\mathcal{C}(\Lambda ,u)$ having that endpoint
is completely inside the loop, whence a contradiction since $C$ and the
other curves of $\mathcal{C}(\Lambda ,u)$ are equivalent up to isometry.

The following properties are true since, if an $\infty $-folding t-curve or
a complete folding t-curve $C$ is the inductive limit of some $n$-folding
t-curves $C_{n}$, then any segment $\left[ x,y\right] $ is contained in $%
\mathrm{F}_{\mathrm{L}}(C)$ (resp. $\mathrm{F}_{\mathrm{R}}(C)$) if and only
if it is contained in $\mathrm{F}_{\mathrm{L}}(C_{n})$ (resp. $\mathrm{F}_{%
\mathrm{R}}(C_{n})$) for $n$ large enough.

If $C$ is an $\infty $-folding t-curve with initial point $x$, then we have $%
\mathrm{F}_{\mathrm{L}}(C)\cap \mathrm{F}_{\mathrm{R}}(C)=\left\{ x\right\} $
and $\mathrm{F}_{\mathrm{L}}(C)$, $\mathrm{F}_{\mathrm{R}}(C)$ are half
curves with endpoint $\emph{x}$.\ If $C$ is a complete folding t-curve, then 
$\mathrm{F}_{\mathrm{L}}(C)$, $\mathrm{F}_{\mathrm{R}}(C)$ are disjoint and
each of them is a complete curve or empty.

For each folding t-curve $C$, we write $\mathrm{F}(C)=\mathrm{F}_{\mathrm{L}%
}(C)\cup \mathrm{F}_{\mathrm{R}}(C)$. We denote by $\mathrm{Dom}(C)$ the
closed part of the plane limited by $\mathrm{F}(C)$ and containing $C$. The
interior of $\mathrm{Dom}(C)$ is connected and $C$ contains the sides of
triangles of $\mathcal{P}$\ which are in $\mathrm{Dom}(C)$\ and not in $%
\mathrm{F}(C)$.

For any $\Lambda ,x$, the sets $\mathrm{Dom}(C)$\ for $C\in \mathcal{C}%
(\Lambda ,x)$ are nonoverlapping and cover the plane.\ If $(C_{i})_{i\in I}$
is a covering of the plane by complete folding t-curves, then the sets $%
\mathrm{Dom}(C_{i})$\ are nonoverlapping and cover the plane, except
possibly the $1$-triangles $(u,v,w)$\ such that $[u,v]$, $[u,w]$, $[v,w]$\
belong to three different curves. We shall see later that, actually, no such
triangle exists.

Now let $n\geq 1$ be an integer and let $C$ be an $n$-folding t-curve
associated to a sequence $T_{\lambda _{1}\cdots \lambda _{n}}$. Consider the
initial points $w,x$ of $C^{\mathrm{I}},C^{\mathrm{M}}$ and the terminal
points $y,z$ of $C^{\mathrm{M}},C^{\mathrm{S}}$.

If $\lambda _{n}=+1$\ (resp. $-1$), then $y\in \mathrm{F}_{\mathrm{L}}(C)$\
and $x\in \mathrm{F}_{\mathrm{R}}(C)$\ (resp. $x\in \mathrm{F}_{\mathrm{L}%
}(C)$\ and $y\in \mathrm{F}_{\mathrm{R}}(C)$). We denote by $\mathrm{F}_{%
\mathrm{LI}}(C)$ the part of $\mathrm{F}_{\mathrm{L}}(C)$\ between $w$ and $%
y $ (resp. $x$), $\mathrm{F}_{\mathrm{LS}}(C)$ the part of $\mathrm{F}_{%
\mathrm{L}}(C)$\ between $y$ (resp. $x$) and $z$, $\mathrm{F}_{\mathrm{RI}%
}(C)$ the part of $\mathrm{F}_{\mathrm{R}}(C)$\ between $w$ and $x$ (resp. $%
y $), $\mathrm{F}_{\mathrm{RS}}(C)$ the part of $\mathrm{F}_{\mathrm{R}}(C)$%
\ between $x$ (resp. $y$) and $z$.

Now suppose $n\geq 2$. If $\lambda _{n}=+1$, then we have $\mathrm{F}_{%
\mathrm{LI}}(C)=\mathrm{F}_{\mathrm{LI}}(C^{\mathrm{I}})\mathrm{F}_{\mathrm{%
LS}}(C^{\mathrm{M}})$, $\mathrm{F}_{\mathrm{LS}}(C)=\mathrm{F}_{\mathrm{L}%
}(C^{\mathrm{S}})$, $\mathrm{F}_{\mathrm{RI}}(C)=\mathrm{F}_{\mathrm{R}}(C^{%
\mathrm{I}})$, $\mathrm{F}_{\mathrm{RS}}(C)=\mathrm{F}_{\mathrm{RI}}(C^{%
\mathrm{M}})\mathrm{F}_{\mathrm{RS}}(C^{\mathrm{S}})$.\ If $\lambda _{n}=-1$%
, then we have $\mathrm{F}_{\mathrm{LI}}(C)=\mathrm{F}_{\mathrm{L}}(C^{%
\mathrm{I}})$, $\mathrm{F}_{\mathrm{LS}}(C)=\mathrm{F}_{\mathrm{LI}}(C^{%
\mathrm{M}})\mathrm{F}_{\mathrm{LS}}(C^{\mathrm{S}})$, $\mathrm{F}_{\mathrm{%
RI}}(C)=\mathrm{F}_{\mathrm{RI}}(C^{\mathrm{I}})\mathrm{F}_{\mathrm{RS}}(C^{%
\mathrm{M}})$, $\mathrm{F}_{\mathrm{RS}}(C)=\mathrm{F}_{\mathrm{R}}(C^{%
\mathrm{S}})$.\bigskip

\noindent \textbf{Lemma 2.8.} Consider $n\in 
%TCIMACRO{\U{2115} }%
%BeginExpansion
\mathbb{N}
%EndExpansion
^{\ast }$ and $\lambda _{1},\ldots ,\lambda _{n}\in \left\{ -1,+1\right\} $.
Let $C$ be an $n$-folding t-curve associated to $T_{\lambda _{1}\cdots
\lambda _{n}}$. Then there exist some sequences $(x_{i})_{0\leq i\leq 2^{n}}$
and $(y_{i})_{0\leq i\leq 2^{n}}$, with $x_{0}=y_{0}$\ initial point of $C$
and $x_{2^{n}}=y_{2^{n}}$\ terminal point of $C$, such that:

\noindent 1) each segment $\left[ x_{i},x_{i+1}\right] $\ or $\left[
y_{i},y_{i+1}\right] $\ is a side of a triangle of $\mathcal{P}$;

\noindent 2) each angle $\widehat{x_{i-1}x_{i}x_{i+1}}$\ or $\widehat{%
y_{i-1}y_{i}y_{i+1}}$\ is equal to $\mp 2\pi /3$;

\noindent 3) $\mathrm{F}_{\mathrm{LI}}(C)=\cup _{1\leq i\leq 2^{n-1}}\left[
x_{i-1},x_{i}\right] $, $\mathrm{F}_{\mathrm{LS}}(C)=\cup _{2^{n-1}+1\leq
i\leq 2^{n}}\left[ x_{i-1},x_{i}\right] $, $\mathrm{F}_{\mathrm{RI}}(C)=\cup
_{1\leq i\leq 2^{n-1}}\left[ y_{i-1},y_{i}\right] $, $\mathrm{F}_{\mathrm{RS}%
}(C)=\cup _{2^{n-1}+1\leq i\leq 2^{n}}\left[ y_{i-1},y_{i}\right] $;

\noindent 4) for $0\leq i\leq 2^{n}$ and $1\leq k\leq n$, each point $%
x_{i},y_{i}$\ belongs to $V_{k}(C)$ if and only if\ $i$ is divisible by $%
2^{k}$.

Now associate to $C$ the sequences $(\alpha _{i})_{1\leq i\leq 2^{n}-1}$ and 
$(\beta _{i})_{1\leq i\leq 2^{n}-1}$ with $\alpha _{i}=+1$ (resp. $-1)$\ if\ 
$\widehat{x_{i-1}x_{i}x_{i+1}}=+2\pi /3$\ (resp. $-2\pi /3$), and $\beta
_{i}=+1$ (resp. $-1)$\ if\ $\widehat{y_{i-1}y_{i}y_{i+1}}=+2\pi /3$\ (resp. $%
-2\pi /3$). Then we have\ $\alpha _{2^{n-1}}=-1$, $\beta _{2^{n-1}}=+1$ and $%
\alpha _{2^{k}+2^{k+1}i}=\beta _{2^{k}+2^{k+1}i}=(-1)^{i}\lambda _{k+2}$ for 
$0\leq k\leq n-2$ and $0\leq i\leq 2^{n-k-1}-1$.\bigskip

\noindent \textbf{Proof.} Lemma 2.8 is clearly true for $n=1$.\ Now suppose
that it is true for an integer $n\geq 1$ and consider an $(n+1)$-folding\
t-curve\ $C$ associated to a sequence $T_{\lambda _{1}\cdots \lambda _{n+1}}$%
. We can assume $\lambda _{n+1}=+1$ since the case $\lambda _{n+1}=-1$ is
similar. We only show the results for $\mathrm{F}_{\mathrm{L}}(C)$ since the
proof for $\mathrm{F}_{\mathrm{R}}(C)$ is similar.

The curves $C^{\mathrm{I}},C^{\mathrm{M}},C^{\mathrm{S}}$\ are associated to 
$T_{\lambda _{1}\cdots \lambda _{n}}$. We have $\mathrm{F}_{\mathrm{LI}}(C)=%
\mathrm{F}_{\mathrm{LI}}(C^{\mathrm{I}})\mathrm{F}_{\mathrm{LS}}(C^{\mathrm{M%
}})$ and $\mathrm{F}_{\mathrm{LS}}(C)=\mathrm{F}_{\mathrm{L}}(C^{\mathrm{S}%
})=\mathrm{F}_{\mathrm{LI}}(C^{\mathrm{S}})\mathrm{F}_{\mathrm{LS}}(C^{%
\mathrm{S}})$. The sequences associated to $\mathrm{F}_{\mathrm{LI}}(C^{%
\mathrm{I}})$\ and $\mathrm{F}_{\mathrm{LI}}(C^{\mathrm{S}})$ are equal. The
sequences associated to $\mathrm{F}_{\mathrm{LS}}(C^{\mathrm{M}})$\ and $%
\mathrm{F}_{\mathrm{LS}}(C^{\mathrm{S}})$ are also equal.\ We turn with a $%
+2\pi /3$ (resp. $-2\pi /3$, $-2\pi /3$) angle when we pass from\ $\mathrm{F}%
_{\mathrm{LI}}(C^{\mathrm{I}})$ to $\mathrm{F}_{\mathrm{LS}}(C^{\mathrm{M}})$
(resp. from\ $\mathrm{F}_{\mathrm{LS}}(C^{\mathrm{M}})$ to $\mathrm{F}_{%
\mathrm{LI}}(C^{\mathrm{S}})$, from\ $\mathrm{F}_{\mathrm{LI}}(C^{\mathrm{S}%
})$ to $\mathrm{F}_{\mathrm{LS}}(C^{\mathrm{S}})$) because $C$ passes twice
(resp. once, once) through the terminal point of $\mathrm{F}_{\mathrm{LI}%
}(C^{\mathrm{I}})$\ (resp. $\mathrm{F}_{\mathrm{LS}}(C^{\mathrm{M}})$, $%
\mathrm{F}_{\mathrm{LI}}(C^{\mathrm{S}})$).\ The results follow from these
facts and from the induction hypothesis applied to $C^{\mathrm{I}},C^{%
\mathrm{M}},C^{\mathrm{S}}$.~~$\blacksquare $\bigskip

\noindent \textbf{Corollary 2.9.} Let $C$ be a complete folding t-curve such
that $\mathrm{F}_{\mathrm{L}}(C)$ is nonempty. Consider the sequence $%
(x_{n})_{n\in 
%TCIMACRO{\U{2124} }%
%BeginExpansion
\mathbb{Z}
%EndExpansion
}$ of endpoints of segments of $\mathrm{F}_{\mathrm{L}}(C)$ and the sequence 
$(\alpha _{n})_{n\in 
%TCIMACRO{\U{2124} }%
%BeginExpansion
\mathbb{Z}
%EndExpansion
}$, with $\alpha _{n}=+1$ if\ $\widehat{x_{n-1}x_{n}x_{n+1}}=+2\pi /3$ and $%
\alpha _{n}=-1$ if\ $\widehat{x_{n-1}x_{n}x_{n+1}}=-2\pi /3$. Then there
exists a sequence $(n_{k})_{k\in 
%TCIMACRO{\U{2115} }%
%BeginExpansion
\mathbb{N}
%EndExpansion
^{\ast }}$\ such that $x_{n_{k}+2^{k}i}\in W_{k}(C)$ and $\alpha
_{n_{k}+2^{k}i}=(-1)^{i}\alpha _{n_{k}}$ for each $k\in 
%TCIMACRO{\U{2115} }%
%BeginExpansion
\mathbb{N}
%EndExpansion
^{\ast }$ and each $i\in 
%TCIMACRO{\U{2124} }%
%BeginExpansion
\mathbb{Z}
%EndExpansion
$. The same property is true for $\mathrm{F}_{\mathrm{R}}(C)$ if it is
nonempty.\bigskip

\noindent \textbf{Remark.} If $C$ is obtained from an $\infty $-folding
t-curve, then there exists a unique $n\in 
%TCIMACRO{\U{2124} }%
%BeginExpansion
\mathbb{Z}
%EndExpansion
-\cup _{k\in 
%TCIMACRO{\U{2115} }%
%BeginExpansion
\mathbb{N}
%EndExpansion
^{\ast }}(n_{k}+2^{k}%
%TCIMACRO{\U{2124} }%
%BeginExpansion
\mathbb{Z}
%EndExpansion
)$; we have $x_{n}\in \cap _{k\in 
%TCIMACRO{\U{2115} }%
%BeginExpansion
\mathbb{N}
%EndExpansion
^{\ast }}V_{k}(C)$. Otherwise, we have $%
%TCIMACRO{\U{2124} }%
%BeginExpansion
\mathbb{Z}
%EndExpansion
=\cup _{k\in 
%TCIMACRO{\U{2115} }%
%BeginExpansion
\mathbb{N}
%EndExpansion
^{\ast }}(n_{k}+2^{k}%
%TCIMACRO{\U{2124} }%
%BeginExpansion
\mathbb{Z}
%EndExpansion
)$.\bigskip

\noindent \textbf{Lemma 2.10.} Any set of disjoint complete folding t-curves
is finite.\bigskip

\noindent \textbf{Proof.} We assume that the sides of triangles of $\mathcal{%
P}$ have length $1$ and we show by induction on $n\in 
%TCIMACRO{\U{2115} }%
%BeginExpansion
\mathbb{N}
%EndExpansion
^{\ast }$\ that $d(x,y)\leq \rho _{n}=\left[ (\sqrt{3})^{n-1}(4-\sqrt{3})-1%
\right] /(\sqrt{3}-1)$ for any points $x,y$ of an\ $n$-folding t-curve. We
have $\rho _{1}=\sqrt{3}$\ and the property is clearly true for $n=1$.

We prove that, if the property is true for an integer $n\geq 1$, then it is
also true for $n+1$. We consider an $(n+1)$-folding t-curve $C$. We
represent $D=\Delta (C)$ in such a way that $C$ and $D$ have the same
endpoints, which gives the length $\sqrt{3}$\ to the segments of $D$.

By the induction hypothesis, for $D$ represented in that way, we have $%
d(x,y)\leq \rho _{n}\sqrt{3}$\ for any $x,y\in D$. Moreover, for each $z\in
C $, there exists $y\in D$\ such that $d(y,z)\leq 1/2$. It follows that, for
any $v,w\in C$, we have $d(v,w)\leq \rho _{n}\sqrt{3}+1=\rho _{n+1}$.

Now we see that, for $x\in 
%TCIMACRO{\U{211d} }%
%BeginExpansion
\mathbb{R}
%EndExpansion
^{2}$ and $k\in 
%TCIMACRO{\U{2115} }%
%BeginExpansion
\mathbb{N}
%EndExpansion
$ large enough, there exists no $r\in 
%TCIMACRO{\U{211d} }%
%BeginExpansion
\mathbb{R}
%EndExpansion
^{+}$ such that $B(x,r)$\ contains points of $k$ disjoint $n$-folding
t-curves for each $n\in 
%TCIMACRO{\U{2115} }%
%BeginExpansion
\mathbb{N}
%EndExpansion
^{\ast }$. This follows since each $n$-folding t-curve contains $3^{n}$
sides of triangles of $\mathcal{P}$, each such curve containing a point of $%
B(x,r)$ is necessarily contained in $B(x,r+\rho _{n})$ and, for $n$ large
compared to $r$, the number of sides of triangles of $\mathcal{P}$ contained
in $B(x,r+\rho _{n})$ is smaller than $3^{n}k$ for an integer $k$ not
depending on $n$ and $r$.~~$\blacksquare $\bigskip

\noindent \textbf{Lemma 2.11.} Let $\mathcal{C}$ be a covering of the plane
by complete folding t-curves. Then the sets $\mathrm{Dom}(C)$ for $C\in 
\mathcal{C}$\ are nonoverlapping and cover the plane, we have $%
V_{k}(C)=V_{k}(D)$ for each $k\in 
%TCIMACRO{\U{2115} }%
%BeginExpansion
\mathbb{N}
%EndExpansion
^{\ast }$\ and any $C,D\in \mathcal{C}$, and all the curves of $\mathcal{C}$
are associated to the same sequence $(\lambda _{n})_{n\in 
%TCIMACRO{\U{2115} }%
%BeginExpansion
\mathbb{N}
%EndExpansion
^{\ast }}$.\bigskip

\noindent \textbf{Proof.} We saw that the sets $\mathrm{Dom}(C)$\ for $C\in 
\mathcal{C}$ are nonoverlapping and cover the plane, except possibly the $1$%
-triangles $(u,v,w)$ such that $[u,v]$, $[u,w]$, $[v,w]$\ belong to three
different curves. Suppose that such a triangle exists and consider the
curves $C,D\in \mathcal{C}$ such that $[u,v]\in C$ and $[u,w]\in D$.

Then we have $\left[ u,v\right] \subset F(C)$ because $w$ cannot belong to $%
C $ and $\left[ u,w\right] \subset F(D)$ because $v$ cannot belong to $D$.
The value of the angle covered by $\mathrm{Dom}(C)$ (resp. $\mathrm{Dom}(D)$%
) in $u$ is $2\pi /3$ if $C$ (resp. $D$) passes once through $u$, and $4\pi
/3$ if $C$ (resp. $D$) passes twice through $u$. If $C$ and $D$ pass only
once through $u$, then another curve $E\in \mathcal{C}$ passes through $u$
and the value of the angle covered by $\mathrm{Dom}(E)$ in $u$ is $2\pi /3$.
In each case, we obtain a contradiction since $\widehat{vuw}$\ and these
angles are all nonoverlapping.

Now we prove the other parts of Lemma 2.11. For each covering $\mathcal{C}$
and each $C\in \mathcal{C}$, we consider the sequence $\Lambda (C)=(\lambda
_{n}(C))_{n\in 
%TCIMACRO{\U{2115} }%
%BeginExpansion
\mathbb{N}
%EndExpansion
^{\ast }}$ associated to $C$. It suffices to show that $V_{1}(C)=V_{1}(D)$
and $\lambda _{1}(C)=\lambda _{1}(D)$ for any $C,D\in \mathcal{C}$, since
these properties imply that $\{\Delta (C)\mid C\in \mathcal{C}\}$\ is also a
covering of the plane by complete folding t-curves. As $C$ is finite by
Lemma 2.10, it is enough to prove these properties for the pairs $(C,D)$\
such that\ $\mathrm{F}(C)\cap \mathrm{F}(D)$\ is unbounded.

For each such pair, we consider some consecutive segments $\left[ x_{0},x_{1}%
\right] ,\ldots ,\left[ x_{8},x_{9}\right] $ of $\mathrm{F}(C)\cap \mathrm{F}%
(D)$, and the integers $\alpha _{1},\ldots ,\alpha _{8}\in \left\{
-1,+1\right\} $ associated to the vertices $x_{1},\ldots ,x_{8}$. According
to Corollary 2.9, we have $\alpha _{1}=-\alpha _{3}=\alpha _{5}=-\alpha _{7}$%
\ or $\alpha _{2}=-\alpha _{4}=\alpha _{6}=-\alpha _{8}$, but the $2$
properties cannot be simultaneously true since the first one implies $\alpha
_{2}=-\alpha _{6}$\ or $\alpha _{4}=-\alpha _{8}$.

Let us suppose for instance that $\alpha _{1}=-\alpha _{3}=\alpha
_{5}=-\alpha _{7}$. Then we have $x_{2}\in V_{1}(C)$ and $x_{2}\in V_{1}(D)$%
, whence $V_{1}(C)=V_{1}(D)$. It follows $\lambda _{1}(C)=\lambda _{1}(D)$
because the $6$ sequences of $3$ segments of curves of $\mathcal{C}$ with
endpoint $x_{2}$ are necessarily equivalent up to rotation.~~$\blacksquare $%
\bigskip

\noindent \textbf{Lemma 2.12.} Here, for each $n\in 
%TCIMACRO{\U{2115} }%
%BeginExpansion
\mathbb{N}
%EndExpansion
$ and each complete folding t-curve $C$, we represent $\Delta ^{n}(C)$ with
segments of length $1$.\ For each $x\in 
%TCIMACRO{\U{211d} }%
%BeginExpansion
\mathbb{R}
%EndExpansion
^{2}$, we denote by $\Delta ^{n}(x)$\ the image of $x$ in this
representation. Then, for each $x\in 
%TCIMACRO{\U{211d} }%
%BeginExpansion
\mathbb{R}
%EndExpansion
^{2}$, any complete folding t-curves $C,D$ and each $n\in 
%TCIMACRO{\U{2115} }%
%BeginExpansion
\mathbb{N}
%EndExpansion
$ large enough, we have $d(\Delta ^{n}(x),\Delta ^{n}(C))<3/4$\ and $%
d(\Delta ^{n}(C),\Delta ^{n}(D))<3/2$.\bigskip

\noindent \textbf{Proof.} For each complete folding t-curve $C$ and each $%
n\in 
%TCIMACRO{\U{2115} }%
%BeginExpansion
\mathbb{N}
%EndExpansion
$, we obtain $\Delta ^{n+1}(C)$ from $\Delta ^{n}(C)$ by replacing each
sequence of $3$ consecutive segments with $1$\ segment, then scaling down
the curve by $\sqrt{3}$. Before the second operation, for each $x\in \Delta
^{n}(C)$, there exists $y\in \Delta ^{n+1}(C)$ such that $d(x,y)\leq 1/2$.

Consequently, for each $x\in 
%TCIMACRO{\U{211d} }%
%BeginExpansion
\mathbb{R}
%EndExpansion
^{2}$ and each $n\in 
%TCIMACRO{\U{2115} }%
%BeginExpansion
\mathbb{N}
%EndExpansion
$, we have

\noindent $d(\Delta ^{n+1}(x),\Delta ^{n+1}(C))\leq (d(\Delta ^{n}(x),\Delta
^{n}(C))+1/2)/\sqrt{3}$, and therefore

\noindent $d(\Delta ^{n+1}(x),\Delta ^{n+1}(C))\leq (5\sqrt{3}/9)d(\Delta
^{n}(x),\Delta ^{n}(C))$ if $d(\Delta ^{n}(x),\Delta ^{n}(C))\geq 3/4$. It
follows $d(\Delta ^{n}(x),\Delta ^{n}(C))<3/4$ for $n$ large enough.

The statement for the curves $C,D$ is an immediate consequence.~~$%
\blacksquare $\bigskip

Now we consider a sequence $X=(x_{n})_{n\in 
%TCIMACRO{\U{2115} }%
%BeginExpansion
\mathbb{N}
%EndExpansion
}\subset U$ with $x_{n+1}\in V_{n}(x_{n})$ for each $n\in 
%TCIMACRO{\U{2115} }%
%BeginExpansion
\mathbb{N}
%EndExpansion
$, a sequence $\Lambda =(\lambda _{n})_{n\in 
%TCIMACRO{\U{2115} }%
%BeginExpansion
\mathbb{N}
%EndExpansion
^{\ast }}\in \left\{ -1,+1\right\} ^{%
%TCIMACRO{\U{2115} }%
%BeginExpansion
\mathbb{N}
%EndExpansion
^{\ast }}$ and the sequences $\Lambda _{n}=(\lambda _{1},\ldots ,\lambda
_{n})$ for $n\in 
%TCIMACRO{\U{2115} }%
%BeginExpansion
\mathbb{N}
%EndExpansion
$. For each $n\in 
%TCIMACRO{\U{2115} }%
%BeginExpansion
\mathbb{N}
%EndExpansion
$, each curve of $\mathcal{C}(\Lambda _{n+1},x_{n+1})$ is obtained by
concatenation of $3$ curves of $\mathcal{C}(\Lambda _{n},x_{n})$. We denote
by $\mathcal{C}(\Lambda ,X)$ the set of inductive limits of curves $C_{n}\in 
\mathcal{C}(\Lambda _{n},x_{n})$.

If there exists $x\in U$ such that $\cap _{n\in 
%TCIMACRO{\U{2115} }%
%BeginExpansion
\mathbb{N}
%EndExpansion
}V_{n}(x_{n})=\left\{ x\right\} $, then $\mathcal{C}(\Lambda ,X)$ contains $%
3 $ half curves starting at $x$ and $3$ reversed such curves ending at $x$.
We denote by $\mathcal{C}^{+}(\Lambda ,X)$ (resp. $\mathcal{C}^{-}(\Lambda
,X)$) the set of curves obtained from $\mathcal{C}(\Lambda ,X)$ by
connecting each terminal segment of reversed half curve with the initial
segment of half curve just at its left (resp. right).

Otherwise, we have $\cap _{n\in 
%TCIMACRO{\U{2115} }%
%BeginExpansion
\mathbb{N}
%EndExpansion
}V_{n}(x_{n})=\emptyset $. Then we write $\mathcal{C}^{0}(\Lambda ,X)=%
\mathcal{C}(\Lambda ,X)$. By Proposition 2.5, each $\mathcal{C}^{\alpha
}(\Lambda ,X)$ is a covering of the plane by complete folding t-curves.

Now, for each $y\in U$ and each $k\in 
%TCIMACRO{\U{2115} }%
%BeginExpansion
\mathbb{N}
%EndExpansion
^{\ast }$, we denote by $H(y,k)$ the regular hexagon of center $y$ such that
each of its sides is the union of $k$ sides of triangles of $\mathcal{P}$,
and $H^{\ast }(y,k)$ its interior.\bigskip

\noindent \textbf{Proposition 2.13.} Consider a sequence $\Lambda =(\lambda
_{k})_{k\in 
%TCIMACRO{\U{2115} }%
%BeginExpansion
\mathbb{N}
%EndExpansion
^{\ast }}\in \left\{ -1,+1\right\} ^{%
%TCIMACRO{\U{2115} }%
%BeginExpansion
\mathbb{N}
%EndExpansion
^{\ast }}$ and two sequences $X=(x_{k})_{k\in 
%TCIMACRO{\U{2115} }%
%BeginExpansion
\mathbb{N}
%EndExpansion
}$, $Y=(y_{k})_{k\in 
%TCIMACRO{\U{2115} }%
%BeginExpansion
\mathbb{N}
%EndExpansion
}$\ such that $x_{k+1}\in V_{k}(x_{k})$ and $y_{k+1}\in V_{k}(y_{k})$ for
each $k\in 
%TCIMACRO{\U{2115} }%
%BeginExpansion
\mathbb{N}
%EndExpansion
$.\ Then, for any $\alpha ,\beta \in \left\{ 0,-,+\right\} $ such that $%
\mathcal{C}^{\alpha }(\Lambda ,X)$\ and $\mathcal{C}^{\beta }(\Lambda ,Y)$\
exist, each $n\in 
%TCIMACRO{\U{2115} }%
%BeginExpansion
\mathbb{N}
%EndExpansion
^{\ast }$ and any $x,y\in U$, there exists $z\in U\cap H(y,5.3^{n})$\ such
that $\mathcal{C}^{\alpha }(\Lambda ,X)\upharpoonright H^{\ast
}(x,3^{n})\cong \mathcal{C}^{\beta }(\Lambda ,Y)\upharpoonright H^{\ast
}(z,3^{n})$.\bigskip

\noindent \textbf{Proof.} There exists $t\in V_{2n+2}(x_{2n+2})$ such that $%
x\in H(t,2.3^{n})$, and therefore $H^{\ast }(x,3^{n})\subset H^{\ast
}(t,3^{n+1})$. For each $u\in V_{2n+2}(y_{2n+2})$, the sets $\mathcal{C}%
^{\alpha }(\Lambda ,X)\upharpoonright H^{\ast }(t,3^{n+1})$ and $\mathcal{C}%
^{\beta }(\Lambda ,Y)\upharpoonright H^{\ast }(u,3^{n+1})$ are isomorphic,
except possibly concerning the way to connect the $6$ segments with endpoint 
$t$ and the $6$ segments with endpoint $u$.

There exist $u,v,w\in V_{2n+2}(y_{2n+2})$ which form a $3^{n+1}$-triangle
containing $y$, which implies $u,v,w\in H(y,3^{n+1})$. One of these points,
say $w$, belongs to $V_{2n+3}(y_{2n+3})$, while the $2$ others belong to $%
W_{2n+3}(y_{2n+3})$. Then the connexions of the $6$ segments of $\mathcal{C}%
^{\beta }(\Lambda ,Y)$ in $u$ and $v$ are different.

Suppose for instance that the connexions of\ $\mathcal{C}^{\alpha }(\Lambda
,X)$ in $t$ are the same as the connexions of\ $\mathcal{C}^{\beta }(\Lambda
,Y)$ in $u$. Then we have\ $\mathcal{C}^{\alpha }(\Lambda ,X)\upharpoonright
H^{\ast }(t,3^{n+1})\cong \mathcal{C}^{\beta }(\Lambda ,Y)\upharpoonright
H^{\ast }(u,3^{n+1})$. Consequently, there exists $z\in U\cap H(u,2.3^{n})$
such that $\mathcal{C}^{\alpha }(\Lambda ,X)\upharpoonright H^{\ast
}(x,3^{n})\cong \mathcal{C}^{\beta }(\Lambda ,Y)\upharpoonright H^{\ast
}(z,3^{n})$. We have $z\in U\cap H(y,5.3^{n})$.~~$\blacksquare $\bigskip

\noindent \textbf{Proof of Theorems 2.1, 2.2, 2.3.} It suffices to prove the
results for oriented curves, since they imply the results for nonoriented
curves.

For the proof of Theorem 2.1, it is enough to consider a curve $C$ with
segments in $E_{i}$. Then it follows from Proposition 2.5 that $C$ belongs
to a covering $\mathcal{C}^{\alpha }(\Lambda ,X)$. Lemma 2.11 implies that $%
\mathcal{C}^{\alpha }(\Lambda ,X)$ is the only covering of the plane by
oriented complete folding t-curves which contains $C$ and satisfies (P). It
satisfies the strong local isomorphism property by Proposition 2.13.
Conversely, by Corollary 2.7, for each covering of the plane by oriented
complete folding t-curves, the local isomorphism property implies (P).

If $\mathcal{C}$ is a covering of the plane by oriented complete folding
t-curves which satisfies (P), then, for any $x,y\in U$, the orientation of
the $6$\ segments with endpoint $x$ and the orientation of the $6$\ segments
with endpoint $y$ are equivalent up to translation. Consequently, if two
such coverings are locally isomorphic, then they have the same orientation.
We see from 3) of Proposition 2.5 that they are associated to the same
sequence $\Lambda $. Conversely, Proposition 2.13 implies that any two such
coverings are locally isomorphic if they have the same orientation and if
they are associated to the same sequence $\Lambda $.~~$\blacksquare $\bigskip

\noindent \textbf{Proof of Theorem 2.4.} Throughout the proof, we use the
notations and the properties introduced just before Lemma 2.8.

First we show that each covering $\mathcal{C}$ contains no more than $3$
curves. By Lemma 2.11, there exist $\alpha \in \left\{ 0,-,+\right\} $ and $%
X=(x_{n})_{n\in 
%TCIMACRO{\U{2115} }%
%BeginExpansion
\mathbb{N}
%EndExpansion
}\subset U$ such that $\mathcal{C}$ is obtained from $\mathcal{C}^{\alpha
}(\Lambda ,X)$ by forgetting its orientation.

If $\cap _{n\in 
%TCIMACRO{\U{2115} }%
%BeginExpansion
\mathbb{N}
%EndExpansion
}V_{n}(x_{n})=\left\{ x\right\} $ with $x\in U$, then $\mathcal{C}%
^{+}(\Lambda ,X)$ (resp. $\mathcal{C}^{-}(\Lambda ,X)$) only contains the $3$
curves $C_{1},C_{2},C_{3}$\ passing through $x$: For each $n\in 
%TCIMACRO{\U{2115} }%
%BeginExpansion
\mathbb{N}
%EndExpansion
$, all the segments\ with points in the interior of $B(\Delta ^{n}(x),1)$\
belong to the $6$\ $1$-folding subcurves of $\Delta ^{n}(C_{1})$, $\Delta
^{n}(C_{2})$, $\Delta ^{n}(C_{3})$ which have $x$ as an endpoint. It follows
that any other curve $D$\ in $\mathcal{C}^{+}(\Lambda ,X)$ (resp. $\mathcal{C%
}^{-}(\Lambda ,X)$) satisfies $d(\Delta ^{n}(x),\Delta ^{n}(D))\geq 1$ for
each $n\in 
%TCIMACRO{\U{2115} }%
%BeginExpansion
\mathbb{N}
%EndExpansion
$, contrary to Lemma 2.12.

If $\cap _{n\in 
%TCIMACRO{\U{2115} }%
%BeginExpansion
\mathbb{N}
%EndExpansion
}V_{n}(x_{n})=\varnothing $ and if $\mathcal{C}^{0}(\Lambda ,X)$ contains
more than $3$ curves, then it contains $4$ curves $C,D,E,F$ with $D$\
separating $C,E$ and $E$ separating $D,F$. For each $n\in 
%TCIMACRO{\U{2115} }%
%BeginExpansion
\mathbb{N}
%EndExpansion
$, $\Delta ^{n}(D)$ is separating $\Delta ^{n}(C)$, $\Delta ^{n}(E)$ and $%
\Delta ^{n}(E)$ is separating $\Delta ^{n}(D)$, $\Delta ^{n}(F)$. It follows 
$d(\Delta ^{n}(C),\Delta ^{n}(F))\geq \sqrt{3}$ since, for any consecutive
segments $r,s,t\in \Delta ^{n}(C)\cup \Delta ^{n}(D)\cup \Delta ^{n}(E)\cup
\Delta ^{n}(F)$, if $r\in \Delta ^{n}(C)$, then $s\in \Delta ^{n}(C)\cup
\Delta ^{n}(D)$, $t\in \Delta ^{n}(C)\cup \Delta ^{n}(D)\cup \Delta ^{n}(E)$%
, and therefore $s,t\notin \Delta ^{n}(F)$. This contradicts Lemma 2.12 for $%
n$ large.

Now we consider the coverings $\mathcal{C}$ associated to $\Lambda $\ which
contain $3$ curves. For each such $\mathcal{C}$, there exists $C,D,E\in 
\mathcal{C}$ distinct such that $F(C)\cap D\neq \varnothing $ and $F(C)\cap
E\neq \varnothing $.

If $F(C)$ consists of $1$ complete curve, then there exists $x\in F(C)$
which belongs to $C,D,E$. We have $x\in \cap _{n\in 
%TCIMACRO{\U{2115} }%
%BeginExpansion
\mathbb{N}
%EndExpansion
}V_{n}(\mathcal{C})$.

Conversely, 3) is realized by any covering $\mathcal{C}$ associated to $%
\Lambda $\ such that there exists $x\in \cap _{n\in 
%TCIMACRO{\U{2115} }%
%BeginExpansion
\mathbb{N}
%EndExpansion
}V_{n}(\mathcal{C})$. Then $\mathcal{C}$\ contains $3$ curves, the
isomorphism class of $\mathcal{C}$ does not depend on $x$ and $\mathcal{C}$
is invariant through a rotation of center $x$ and angle $2\pi /3$. We obtain
a representative of the other isomorphism class of coverings of that type by
changing the connections in $x$.

If $F(C)$ consists of $2$ complete curves, then we can orient $C,D,E$ in
such a way that $F_{L}(C)=F_{L}(D)$ and $F_{R}(C)=F_{R}(E)$. We consider $2$
segments $s\subset F_{L}(C)$, $t\subset F_{R}(C)$ of length $1$, and $2$
sequences $(C_{n})_{n\in 
%TCIMACRO{\U{2115} }%
%BeginExpansion
\mathbb{N}
%EndExpansion
}$, $(P_{n})_{n\in 
%TCIMACRO{\U{2115} }%
%BeginExpansion
\mathbb{N}
%EndExpansion
}$ such that $C=\cup _{n\in 
%TCIMACRO{\U{2115} }%
%BeginExpansion
\mathbb{N}
%EndExpansion
}C_{n}$ and such that, for each $n\in 
%TCIMACRO{\U{2115} }%
%BeginExpansion
\mathbb{N}
%EndExpansion
$, $C_{n}$\ is an $n$-folding t-curve, $P_{n}\in \left\{ \mathrm{I},\mathrm{M%
},\mathrm{S}\right\} $\ and $C_{n}=C_{n+1}^{P_{n}}$. There exists an integer 
$h$ such that $s\subset \mathrm{F}_{\mathrm{L}}(C_{n})$\ and $t\subset 
\mathrm{F}_{\mathrm{R}}(C_{n})$\ for each $n\geq h$. We consider the
sequences $(\alpha _{n})_{n\geq h}$ and $(\beta _{n})_{n\geq h}$\ such that $%
\alpha _{n},\beta _{n}\in \left\{ \mathrm{I},\mathrm{S}\right\} $, $s\subset 
\mathrm{F}_{\mathrm{L}\alpha _{n}}(C_{n})$\ and $t\subset \mathrm{F}_{%
\mathrm{R}\beta _{n}}(C_{n})$\ for each $n$.

For $n\geq h$, the properties $s\subset \mathrm{F}_{\mathrm{L}\alpha
_{n}}(C_{n})\cap \mathrm{F}_{\mathrm{L}\alpha _{n+1}}(C_{n+1})$\ and $%
t\subset \mathrm{F}_{\mathrm{R}\beta _{n}}(C_{n})\cap \mathrm{F}_{\mathrm{R}%
\beta _{n+1}}(C_{n+1})$\ respectively imply $\mathrm{F}_{\mathrm{L}\alpha
_{n}}(C_{n+1}^{P_{n}})=\mathrm{F}_{\mathrm{L}\alpha _{n}}(C_{n})\subset 
\mathrm{F}_{\mathrm{L}\alpha _{n+1}}(C_{n+1})$\ and $\mathrm{F}_{\mathrm{R}%
\beta _{n}}(C_{n+1}^{P_{n}})=\mathrm{F}_{\mathrm{R}\beta _{n}}(C_{n})\subset 
\mathrm{F}_{\mathrm{R}\beta _{n+1}}(C_{n+1})$. Considering $C_{n+1}$, we see
that one of the $4$\ cases below is realized:

\noindent $P_{n}=\mathrm{I}$, $(\alpha _{n},\beta _{n})\in \left\{ (\mathrm{I%
},\mathrm{I}),(\mathrm{I},\mathrm{S}),(\mathrm{S},\mathrm{I})\right\} $,\ $%
(\alpha _{n+1},\beta _{n+1})=(\mathrm{I},\mathrm{I})$;

\noindent $P_{n}=\mathrm{S}$, $(\alpha _{n},\beta _{n})\in \left\{ (\mathrm{I%
},\mathrm{S}),(\mathrm{S},\mathrm{I}),(\mathrm{S},\mathrm{S})\right\} $,\ $%
(\alpha _{n+1},\beta _{n+1})=(\mathrm{S},\mathrm{S})$;

\noindent $P_{n}=\mathrm{M}$, $(\alpha _{n},\beta _{n})=(\mathrm{S},\mathrm{I%
})$,\ $(\alpha _{n+1},\beta _{n+1})=(\mathrm{I},\mathrm{S})$, $\lambda
_{n+1}=+1$;

\noindent $P_{n}=\mathrm{M}$, $(\alpha _{n},\beta _{n})=(\mathrm{I},\mathrm{S%
})$,\ $(\alpha _{n+1},\beta _{n+1})=(\mathrm{S},\mathrm{I})$, $\lambda
_{n+1}=-1$.

It follows that there exists $m\geq h$ such that one of the properties below
is true for each $n\geq m$:

\noindent a) $P_{n}=\mathrm{I}$ and $(\alpha _{n},\beta _{n})=(\mathrm{I},%
\mathrm{I})$;

\noindent b) $P_{n}=\mathrm{S}$ and $(\alpha _{n},\beta _{n})=(\mathrm{S},%
\mathrm{S})$;

\noindent c) $P_{n}=\mathrm{M}$; $(\alpha _{n},\beta _{n})=(\mathrm{I},%
\mathrm{S})$ and $\lambda _{n}=+1$ for $n$ even; $(\alpha _{n},\beta _{n})=(%
\mathrm{S},\mathrm{I})$ and $\lambda _{n}=-1$ for $n$ odd;

\noindent d) $P_{n}=\mathrm{M}$; $(\alpha _{n},\beta _{n})=(\mathrm{I},%
\mathrm{S})$ and $\lambda _{n}=+1$ for $n$ odd; $(\alpha _{n},\beta _{n})=(%
\mathrm{S},\mathrm{I})$ and $\lambda _{n}=-1$ for $n$ even.

The cases a) and b) imply $\cup _{n\in 
%TCIMACRO{\U{2115} }%
%BeginExpansion
\mathbb{N}
%EndExpansion
}C_{n}\neq C$, contrary to our hypotheses.

On the other hand, each of the cases c), d), for each of the values of $%
\Lambda $ which realize it, gives a covering of the plane by complete
folding t-curves, with $1$ curve separating $2$ others. This covering is
determined by $\Lambda $ modulo a positive isometry. We obtain
representatives of the $2$ other isomorphism classes of coverings by
applying a rotation of angle $\mp 2\pi /3$.

It remains to be proved that there exist $2^{\omega }$ isomorphism classes
of coverings of the plane by $1$ (resp. $2$) complete folding t-curves
associated to $\Lambda $. We consider the complete folding t-curves $C$ for
which there exists a pair $((C_{n})_{n\in 
%TCIMACRO{\U{2115} }%
%BeginExpansion
\mathbb{N}
%EndExpansion
},(P_{n})_{n\in 
%TCIMACRO{\U{2115} }%
%BeginExpansion
\mathbb{N}
%EndExpansion
})$ such that $C=\cup _{n\in 
%TCIMACRO{\U{2115} }%
%BeginExpansion
\mathbb{N}
%EndExpansion
}C_{n}$ and such that, for each $n\in 
%TCIMACRO{\U{2115} }%
%BeginExpansion
\mathbb{N}
%EndExpansion
$, $C_{n}$\ is an $n$-folding t-curve, $P_{n}\in \left\{ \mathrm{I},\mathrm{M%
},\mathrm{S}\right\} $\ and $C_{n}=C_{n+1}^{P_{n}}$. Then $C$ is determined
modulo a positive isometry by $(P_{n})_{n\in 
%TCIMACRO{\U{2115} }%
%BeginExpansion
\mathbb{N}
%EndExpansion
}$ and there are countably many choices of $(P_{n})_{n\in 
%TCIMACRO{\U{2115} }%
%BeginExpansion
\mathbb{N}
%EndExpansion
}$ which give $C$, since any $2$ such sequences only differ by a finite
number of terms. Consequently, it suffices to show that $2^{\omega }$
choices of $(P_{n})_{n\in 
%TCIMACRO{\U{2115} }%
%BeginExpansion
\mathbb{N}
%EndExpansion
}$ give coverings by $1$ (resp. $2$) complete folding t-curves.

In order to obtain a curve $C$ which covers the plane, it suffices to have $%
\mathrm{F}_{\mathrm{L}}(C_{5n})\cap \mathrm{F}_{\mathrm{L}%
}(C_{5n+2})=\varnothing $ and $\mathrm{F}_{\mathrm{R}}(C_{5n+2})\cap \mathrm{%
F}_{\mathrm{R}}(C_{5n+4})=\varnothing $ for each $n\in 
%TCIMACRO{\U{2115} }%
%BeginExpansion
\mathbb{N}
%EndExpansion
$. The first property is realized for $P_{5n}=\mathrm{S}$, $P_{5n+1}=\mathrm{%
I}$ if $\lambda _{5n+2}=+1$, and $P_{5n}=\mathrm{I}$, $P_{5n+1}=\mathrm{S}$
if $\lambda _{5n+2}=$ $-1$. The second property is realized for $P_{5n+2}=%
\mathrm{I}$, $P_{5n+3}=\mathrm{S}$ if $\lambda _{5n+4}=+1$, and $P_{5n+2}=%
\mathrm{S}$, $P_{5n+3}=\mathrm{I}$ if $\lambda _{5n+4}=-1$. For each $n\in 
%TCIMACRO{\U{2115} }%
%BeginExpansion
\mathbb{N}
%EndExpansion
$, $P_{5n+4}$ can be chosen arbitrarily. It follows that $2^{\omega }$
choices give coverings by $1$ curve.

Now we prove that $2^{\omega }$ choices give coverings by $2$ curves. As
only countably many choices give coverings by $3$ curves, it suffices to
show that $2^{\omega }$ choices give a curve $C$ with $\mathrm{F}_{\mathrm{L}%
}(C)\neq \varnothing $.

For each $n\in 
%TCIMACRO{\U{2115} }%
%BeginExpansion
\mathbb{N}
%EndExpansion
^{\ast }$, we consider the sequences $(P_{1},\ldots ,P_{n-1})$, $%
(C_{1},\ldots ,C_{n})$, $(\alpha _{1},\ldots ,\alpha _{n})$, with $C_{i}$ an 
$i$-folding t-curve for $1\leq i\leq n$, $C_{i}=C_{i+1}^{P_{i}}$ and $%
\mathrm{F}_{\mathrm{L}\alpha _{i}}(C_{i})\subset \mathrm{F}_{\mathrm{L}%
\alpha _{i+1}}(C_{i+1})$\ for $1\leq i\leq n-1$. We observe that, for any
such sequences, there are $2$ different ways to choose $P_{n}$ so that the $%
(n+1)$-folding t-curve $C_{n+1}$\ defined by $C_{n}=C_{n+1}^{P_{n}}$
satisfies $\mathrm{F}_{\mathrm{L}\alpha _{n}}(C_{n})\subset \mathrm{F}_{%
\mathrm{L}}(C_{n+1})$, and therefore $\mathrm{F}_{\mathrm{L}\alpha
_{n}}(C_{n})\subset \mathrm{F}_{\mathrm{L}\alpha _{n+1}}(C_{n+1})$\ for the
appropriate $\alpha _{n+1}$.~~$\blacksquare $\bigskip

\noindent \textbf{Remark.} Each covering of the plane by $3$ nonoriented
complete folding t-curves is invariant through a central symmetry.\bigskip

\textbf{3. Peano-Gosper curves}\bigskip

In the present section, we obtain\ a set $\Omega $ of coverings of $%
%TCIMACRO{\U{211d} }%
%BeginExpansion
\mathbb{R}
%EndExpansion
^{2}$ by sets of disjoint self-avoiding nonoriented curves, generalizing the
Peano-Gosper curves, such that:

\noindent 1) each $\mathcal{C}\in \Omega $ satisfies the strong local
isomorphism property; any set of curves locally isomorphic to $\mathcal{C}$
belongs to $\Omega $;

\noindent 2) $\Omega $ is the union of $2^{\omega }$ equivalence classes for
the relation ``$\mathcal{C}$ locally\ isomorphic to $\mathcal{D}$''; each of
them contains $2^{\omega }$\ (resp. $2^{\omega }$, $1$, $0$)\ isomorphism
classes of coverings by $1$ (resp. $2$, $3$, $\geq 4$) curves.

\noindent Each $\mathcal{C}\in \Omega $ gives exactly $2$ coverings of $%
%TCIMACRO{\U{211d} }%
%BeginExpansion
\mathbb{R}
%EndExpansion
^{2}$ by sets of oriented curves which satisfy the local isomorphism
property. They have opposite orientations.

For each regular tiling $\mathcal{P}$ of $%
%TCIMACRO{\U{211d} }%
%BeginExpansion
\mathbb{R}
%EndExpansion
^{2}$ by hexagons and each center $x$ of a tile of $\mathcal{P}$, we
construct some tilings $\mathcal{P}_{x\lambda _{1}\cdots \lambda _{n}}$ of $%
%TCIMACRO{\U{211d} }%
%BeginExpansion
\mathbb{R}
%EndExpansion
^{2}$ by isomorphic tiles such that, for each $n\in 
%TCIMACRO{\U{2115} }%
%BeginExpansion
\mathbb{N}
%EndExpansion
$\ and\ any\ $\lambda _{1},\ldots ,\lambda _{n}\in \left\{ -,+\right\} $, $x$
is the center of a tile of $\mathcal{P}_{x\lambda _{1}\cdots \lambda _{n}}$
and, for $n\geq 1$, each tile of $\mathcal{P}_{x\lambda _{1}\cdots \lambda
_{n}}$ is the union of $7$ nonoverlapping tiles of $\mathcal{P}_{x\lambda
_{1}\cdots \lambda _{n-1}}$ with one of them surrounded by the $6$ others.

\begin{center}
\includegraphics[scale=0.80]{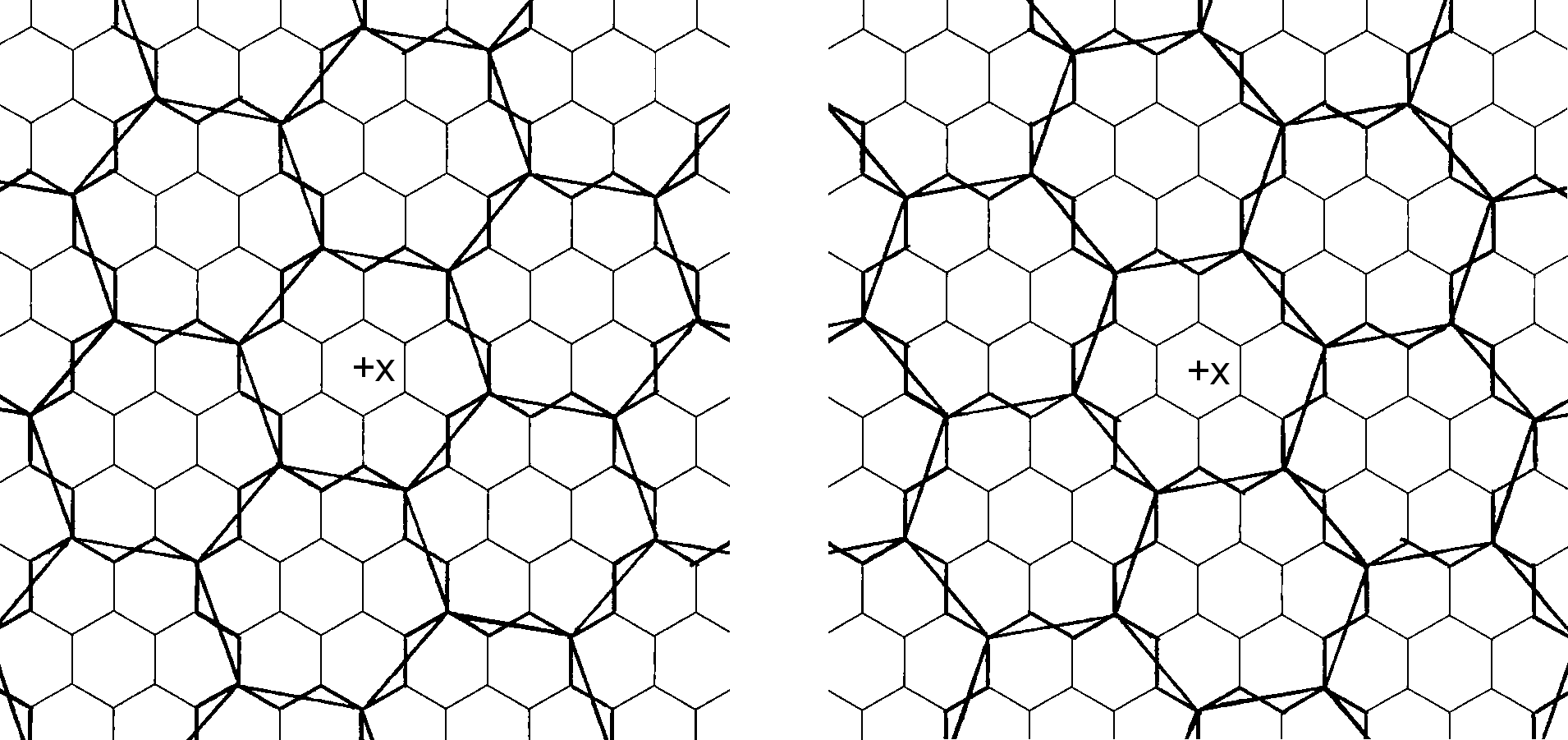}

\medskip Figure 3.1
\end{center}

We write $\mathcal{P}_{x}=\mathcal{P}$. We consider the tilings $\mathcal{P}%
_{x+}$ and $\mathcal{P}_{x-}$ respctively given by the first and the second
part of Figure 3.1. The points of $%
%TCIMACRO{\U{211d} }%
%BeginExpansion
\mathbb{R}
%EndExpansion
^{2}$ which are common to $3$\ tiles of $\mathcal{P}_{x+}$\ (resp. $\mathcal{%
P}_{x-}$) determine a regular tiling by hexagons $\mathcal{Q}_{x+}$\ (resp. $%
\mathcal{Q}_{x-}$).

We denote by $\Delta _{x+}$ (resp. $\Delta _{x-}$) the bijection which
associates to each tile of $\mathcal{P}_{x+}$\ (resp. $\mathcal{P}_{x-}$)
the tile of $\mathcal{Q}_{x+}$\ (resp. $\mathcal{Q}_{x-}$) with the same
center. We see from Figure 3.1 that, for each tile $Q$\ of $\mathcal{Q}_{x+}$%
\ (resp. $\mathcal{Q}_{x-}$), $\Delta _{x+}^{-1}(Q)$ (resp.\ $\Delta
_{x-}^{-1}(Q)$) is obtained by replacing each side $S$ of $Q$ with three
consecutive sides of tiles of $\mathcal{P}$. The first side is at the right
(resp. left) of $S$ and the third side is at its left (resp. right) for the
direction given by anticlockwise rotation around the center of $Q$.

For $n\geq 2$\ and\ $\lambda _{1},\ldots ,\lambda _{n}\in \left\{
-,+\right\} $ we write $\mathcal{P}_{x\lambda _{1}\cdots \lambda
_{n}}=\Delta _{x\lambda _{1}}^{-1}((\mathcal{Q}_{x\lambda _{1}})_{x\lambda
_{2}\cdots \lambda _{n}})$. Each\ $\mathcal{P}_{x\lambda _{1}\cdots \lambda
_{n}}$ consists of isomorphic tiles with connected interior and connected
exterior. It is regular and invariant through a rotation of center $x$ and
angle $\pi /3$.

In each\ $\mathcal{P}_{x\lambda _{1}\cdots \lambda _{n}}$, each tile $P$ has
nonempty frontiers with $6$ others. We call \emph{sides} of $P$ these
frontiers, which are unions of $3^{n}$ sides of tiles of $\mathcal{P}$. The 
\emph{vertices} of $P$ are the endpoints of its sides.

For each $n\in 
%TCIMACRO{\U{2115} }%
%BeginExpansion
\mathbb{N}
%EndExpansion
$, any\ $\lambda _{1},\ldots ,\lambda _{n}\in \left\{ -,+\right\} $\ and any
centers $x,y$ of tiles of $\mathcal{P}$,\ each tile of $\mathcal{P}%
_{x\lambda _{1}\cdots \lambda _{n}-}$ and each tile of $\mathcal{P}%
_{y\lambda _{1}\cdots \lambda _{n}+}$ are isomorphic.

For any integers $n\geq m\geq 1$, we define some \emph{derivations}$\ \Delta
_{x\lambda _{1}\cdots \lambda _{m}}$\ on the sets $\mathcal{P}_{x\lambda
_{1}\cdots \lambda _{n}}$. We use the definitions of $\Delta _{x+}$ and $%
\Delta _{x-}$ given above and we write $\Delta _{x\lambda _{1}\cdots \lambda
_{m}}(P)=\Delta _{x\lambda _{m}}(\cdots (\Delta _{x\lambda _{1}}(P))\cdots )$%
\ for any $\lambda _{1},\ldots ,\lambda _{n}\in \left\{ -,+\right\} $ and
each $P\in \mathcal{P}_{x\lambda _{1}\cdots \lambda _{n}}$. We have

\noindent $\Delta _{x\lambda _{1}\cdots \lambda _{m}}(\mathcal{P}_{x\lambda
_{1}\cdots \lambda _{n}})=(\Delta _{x\lambda _{1}\cdots \lambda _{m}}(%
\mathcal{P}_{x\lambda _{1}\cdots \lambda _{m}}))_{x\lambda _{m+1}\cdots
\lambda _{n}}$.\bigskip

\noindent \textbf{Remark}. For each $n\in 
%TCIMACRO{\U{2115} }%
%BeginExpansion
\mathbb{N}
%EndExpansion
$, denote by $P_{n}$ the tile of center $x$ in $\mathcal{P}_{x\lambda
_{1}\cdots \lambda _{n}}$ with $\lambda _{1}=\cdots =\lambda _{n}=+$, and $%
Q_{n}$ the tile with the same vertices as $P_{0}$ which is the image of $%
P_{n}$ under a similarity. Then the limit of the tiles $Q_{n}$ is the
Peano-Gosper island considered in [3, p. 46]. It is the union of $7$
isomorphic nonoverlapping tiles which are similar to it, with one of them
surrounded by the $6$ others.\bigskip

Now, for each $\Lambda =(\lambda _{n})_{n\in 
%TCIMACRO{\U{2115} }%
%BeginExpansion
\mathbb{N}
%EndExpansion
^{\ast }}\in \left\{ -,+\right\} ^{%
%TCIMACRO{\U{2115} }%
%BeginExpansion
\mathbb{N}
%EndExpansion
^{\ast }}$, we consider the sequences $X=(x_{n})_{n\in 
%TCIMACRO{\U{2115} }%
%BeginExpansion
\mathbb{N}
%EndExpansion
}$\ of centers of tiles of $\mathcal{P}$ with $x_{n+1}=x_{n}$ or $%
x_{n},x_{n+1}$ centers of adjacent tiles of $\mathcal{P}_{x_{n}\lambda
_{1}\cdots \lambda _{n}}$ for each $n\in 
%TCIMACRO{\U{2115} }%
%BeginExpansion
\mathbb{N}
%EndExpansion
$. For each such $X$ and each $n\in 
%TCIMACRO{\U{2115} }%
%BeginExpansion
\mathbb{N}
%EndExpansion
$, we write $\mathcal{P}_{X\Lambda }^{n}=\mathcal{P}_{x_{n}\lambda
_{1}\cdots \lambda _{n}}$ and we denote by $P_{X\Lambda }^{n}$\ the tile of $%
\mathcal{P}_{X\Lambda }^{n}$ which contains $x_{n}$. For $n\in 
%TCIMACRO{\U{2115} }%
%BeginExpansion
\mathbb{N}
%EndExpansion
$, we have $P_{X\Lambda }^{n}\subset P_{X\Lambda }^{n+1}$ and each tile of $%
\mathcal{P}_{X\Lambda }^{n+1}$ is the union of $7$ nonoverlapping tiles of $%
\mathcal{P}_{X\Lambda }^{n}$. We write $\mathcal{P}_{X\Lambda }=\cup _{n\in 
%TCIMACRO{\U{2115} }%
%BeginExpansion
\mathbb{N}
%EndExpansion
}\mathcal{P}_{X\Lambda }^{n}$.

For each center $u$ of a tile of $\mathcal{P}$,\ there exists a unique
sequence $S_{u}=(u_{n})_{n\in 
%TCIMACRO{\U{2115} }%
%BeginExpansion
\mathbb{N}
%EndExpansion
}$ with the properties above\ such that $u_{0}=u$\ and $\mathcal{P}%
_{S_{u}\Lambda }=\mathcal{P}_{X\Lambda }$. We write $R_{uX\Lambda }=\cup
_{n\in 
%TCIMACRO{\U{2115} }%
%BeginExpansion
\mathbb{N}
%EndExpansion
}P_{S_{u}\Lambda }^{n}$.

We call a \emph{region}\ of $\mathcal{P}_{X\Lambda }$ any minimal nonempty
closed subset $R\subset 
%TCIMACRO{\U{211d} }%
%BeginExpansion
\mathbb{R}
%EndExpansion
^{2}$ such that any tile\ of $\mathcal{P}_{X\Lambda }$ is completely inside $%
R$\ or completely outside. Each $R_{uX\Lambda }$ is a region of $\mathcal{P}%
_{X\Lambda }$ and each region of $\mathcal{P}_{X\Lambda }$\ is obtained in
that way.

For each $n\in 
%TCIMACRO{\U{2115} }%
%BeginExpansion
\mathbb{N}
%EndExpansion
$ and any centers $u,v$ of tiles of $\mathcal{P}$, $P_{S_{u}\Lambda }^{n}$\
and $P_{S_{v}\Lambda }^{n}$\ are disjoint, or they have one common side, or
they are equal. If the third possibility is realized for some $n$, then we
have $R_{uX\Lambda }=R_{vX\Lambda }$. Otherwise, the second possibility is
realized for $n$ large enough and $R_{uX\Lambda },R_{vX\Lambda }$\ are
adjacent.

For any centers $u,v,w$ of tiles of $\mathcal{P}$, if\ $P_{S_{u}\Lambda
}^{n},P_{S_{v}\Lambda }^{n},P_{S_{w}\Lambda }^{n}$\ are distinct for each $%
n\in 
%TCIMACRO{\U{2115} }%
%BeginExpansion
\mathbb{N}
%EndExpansion
$, then, for $n$\ large enough, they have $1$ common point and any $2$ of
them have a common side containing that point. Consequently, $R_{uX\Lambda
},R_{vX\Lambda },R_{wX\Lambda }$\ are the three regions of $\mathcal{P}%
_{X\Lambda }$.

It follows that one of the three following properties is true for each $%
\mathcal{P}_{X\Lambda }$:

\noindent 1) There is $1$ region.

\noindent 2) There are $2$ regions; their frontier is a complete curve which
consists of sides of hexagons of $\mathcal{P}$.

\noindent 3) There are $3$ regions with $1$ common point which is a vertex
of tiles of $\mathcal{P}_{X\Lambda }^{n}$ for each $n\in 
%TCIMACRO{\U{2115} }%
%BeginExpansion
\mathbb{N}
%EndExpansion
$; the frontier of $2$ regions is a half curve starting from that point
which consists of sides of hexagons of $\mathcal{P}$.\bigskip

\noindent \textbf{Proposition 3.1.} For each $\Lambda $, the isomorphism\
classes of sets $\mathcal{P}_{X\Lambda }$ are countable; there exist $%
2^{\omega }$ classes of sets with $1$\ region, $2^{\omega }$ classes of sets
with $2$\ regions, and $2$ classes of sets with $3$\ regions, obtained from
each other by a rotation of angle $\pi /3$.\bigskip

\noindent \textbf{Proof.} For each sequence $X=(x_{n})_{n\in 
%TCIMACRO{\U{2115} }%
%BeginExpansion
\mathbb{N}
%EndExpansion
}$\ such that $\mathcal{P}_{X\Lambda }$ exists, we write $R_{X\Lambda
}=R_{x_{0}X\Lambda }$. For any such sequences $X,Y$, we write $X\sim _{1}Y$
if $\mathcal{P}_{X\Lambda }\cong \mathcal{P}_{Y\Lambda }$, $X\sim _{2}Y$ if $%
\mathcal{P}_{X\Lambda }=\mathcal{P}_{Y\Lambda }$, and $X\sim _{3}Y$ if $%
\mathcal{P}_{X\Lambda }=\mathcal{P}_{Y\Lambda }$ and $R_{X\Lambda
}=R_{Y\Lambda }$.

Each $\sim _{3}$-class is countable since any two sequences in such a class
are ultimately equal. As each $\mathcal{P}_{X\Lambda }$\ only has finitely
many regions, each $\sim _{2}$-class is also countable.

For any sequences $X,Y$\ such that $X\sim _{1}Y$, there exist a sequence $Z$
and a translation $\tau $ such that $X\sim _{2}Z$ and $\tau (Z)=Y$; we have $%
\tau (\mathcal{P})=\mathcal{P}$. As only countably many translations satisfy
that property, each $\sim _{1}$-class is the union of countably many $\sim
_{2}$-classes, and therefore countable.

Now we prove the second part of the proposition.

First we consider a vertex $w$ of a tile of $\mathcal{P}$ and the sequences $%
Y$ such that $w$ is a vertex of $P_{Y\Lambda }^{n}$ for each $n\in 
%TCIMACRO{\U{2115} }%
%BeginExpansion
\mathbb{N}
%EndExpansion
$. This property is true for exactly $6$ sequences $Y_{1},\ldots ,Y_{6}$
with $R_{Y_{1}\Lambda },R_{Y_{2}\Lambda },R_{Y_{3}\Lambda }$ distinct
regions of $\mathcal{P}_{Y_{1}\Lambda }=\mathcal{P}_{Y_{2}\Lambda }=\mathcal{%
P}_{Y_{3}\Lambda }$, $R_{Y_{1}\Lambda }\cap R_{Y_{2}\Lambda }\cap
R_{Y_{3}\Lambda }=\left\{ w\right\} $ and $Y_{4},Y_{5},Y_{6}$ obtained from $%
Y_{1},Y_{2},Y_{3}$ by a rotation of center $w$ and angle $\pi /3$.

The statement concerning the classes of sets with $3$ regions follows since,
for each sequence $X$, if $\mathcal{P}_{X\Lambda }$ has $3$ regions with the
common point $w$, then $X$ is ultimately equal to some $Y_{i}$.

Now, it suffices to prove that there exist $2^{\omega }$ sequences $%
X=(x_{n})_{n\in 
%TCIMACRO{\U{2115} }%
%BeginExpansion
\mathbb{N}
%EndExpansion
}$ such that $\mathcal{P}_{X\Lambda }$ has $1$ region and $2^{\omega }$
sequences with $2$ or $3$ regions.

In order to obtain a sequence $X$ with $1$ region, it suffices to choose
successively the elements $x_{k}$ so that $P_{X\Lambda }^{k-1}$ is contained
in the interior of $P_{X\Lambda }^{k}$ for each $k$\ even. Consequently,
there are $2^{\omega }$ such sequences.

For each side $S$ of a tile of $\mathcal{P}$, in order to obtain a sequence $%
X$ with $2$ or $3$ regions, it suffices to choose successively the elements $%
x_{k}$ so that $S$ is contained in a side of $P_{X\Lambda }^{k}$. There are $%
3$ possible choices\ for each $k\geq 1$, and therefore $2^{\omega }$ such
sequences.~~$\blacksquare $\bigskip

Now we define the curves associated to $\mathcal{P}$. The set $V$ of
vertices of hexagons of $\mathcal{P}$ is the disjoint union of $2$ subsets\
which contain no consecutive vertices of an hexagon. We choose one of them
and we denote it by $W$.

An \emph{oriented bounded curve} (resp. \emph{half curve}, \emph{complete
curve}) is a sequence $C=(A_{k})_{0\leq k\leq n}$\ (resp. $(A_{k})_{k\in 
%TCIMACRO{\U{2115} }%
%BeginExpansion
\mathbb{N}
%EndExpansion
}$, $(A_{k})_{k\in 
%TCIMACRO{\U{2124} }%
%BeginExpansion
\mathbb{Z}
%EndExpansion
}$) of consecutive oriented segments, each of them joining $2$\ vertices of
an hexagon of $\mathcal{P}$ which belong to $W$, such that each hexagon
contains at most $1$\ segment\ and each endpoint of $C$, if it exists, only
belongs to $1$ segment. With this definition, $C$ is self-avoiding.

For each oriented bounded curve $(A_{0},\ldots ,A_{n})$, we consider the
sequence\ $(a_{1},\ldots ,a_{n})$\ with each $a_{k}$ equal to $+2,+1,0,-1,-2$%
\ if the angle between $A_{k-1}$ and $A_{k}$ is $+2\pi /3,+\pi /3,0,-\pi
/3,-2\pi /3$. For each sequence\ $S=(a_{1},\ldots ,a_{n})$, we write $%
\overline{S}=(-a_{n},\ldots ,-a_{1})$. We have $\overline{\overline{S}}=S$.
The curves associated to $\overline{S}$\ are obtained by changing the
orientation of the curves associated to $S$.

For each $n\in 
%TCIMACRO{\U{2115} }%
%BeginExpansion
\mathbb{N}
%EndExpansion
$, any\ $\lambda _{1},\ldots ,\lambda _{n}\in \left\{ -,+\right\} $, each
center $w$ of a tile of $\mathcal{P}$ and each $P\in \mathcal{P}_{w\lambda
_{1}\cdots \lambda _{n}}$, we say that a curve $C$ is a \emph{covering} of $%
P $ if it is contained in $P$ and if, for each $m\in \left\{ 0,\ldots
,n\right\} $ and each $Q\in \mathcal{P}_{w\lambda _{1}\cdots \lambda _{m}}$
contained in $P$, the segments of $C$ contained in $Q$ form a subcurve whose
endpoints are\ vertices of $Q$.

For any integers $n\geq m\geq 1$ and each covering $C$ of $P\in \mathcal{P}%
_{w\lambda _{1}\cdots \lambda _{n}}$, we obtain a covering $\Delta
_{w\lambda _{1}\cdots \lambda _{m}}(C)$\ of $\Delta _{w\lambda _{1}\cdots
\lambda _{m}}(P)$ by replacing each $C\upharpoonright Q$ for $Q\in \mathcal{P%
}_{w\lambda _{1}\cdots \lambda _{m}}$\ contained in $P$ with the segment
from its initial point to its terminal point. We have $\Delta _{w\lambda
_{1}\cdots \lambda _{m}}(C)=\Delta _{w\lambda _{m}}(\cdots (\Delta
_{w\lambda _{1}}(C))\cdots )$.

Now we define by induction on $n\in 
%TCIMACRO{\U{2115} }%
%BeginExpansion
\mathbb{N}
%EndExpansion
^{\ast }$ some sequences

\noindent $S_{\lambda _{1}\cdots \lambda _{n}}\in \left\{
-2,-1,0,+1,+2\right\} ^{7^{n}-1}$\ for $\lambda _{1},\ldots ,\lambda _{n}\in
\left\{ -,+\right\} $. We write

\noindent $S_{+}=S=(+1,+2,-1,-2,0,-1)$, $S_{-}=-S$,

\noindent $S_{++}=(S,+1,\overline{S},+1,\overline{S},-1,S,-1,S,+1,S,-1,%
\overline{S})$, $S_{--}=-S_{++}$,

\noindent $S_{+-}=(\overline{S},-1,S,-1,S,+1,\overline{S},+1,\overline{S},-1,%
\overline{S},+1,S)$, $S_{-+}=-S_{+-}$.

\noindent For $n\geq 3$, $S_{\lambda _{1}\cdots \lambda _{n}}$ is obtained
from $S_{\lambda _{2}\cdots \lambda _{n}}=(a_{1},\ldots ,a_{7^{n-1}-1})$ by
replacing each subsequence $(a_{7k+1},\ldots ,a_{7k+6})$ equal to $%
S_{\lambda _{2}}$\ (resp. $\overline{S_{\lambda _{2}}}$) with $S_{\lambda
_{1}\lambda _{2}}$\ (resp. $\overline{S_{\lambda _{1}\lambda _{2}}}$).

The first (resp. second) part of Figure 3.2 below shows a covering of a tile
of $\mathcal{P}_{w+}$ (resp. $\mathcal{P}_{w-}$) by a curve associated to $%
S_{+}$ (resp. $S_{-}$).

\begin{center}
\includegraphics[scale=0.75]{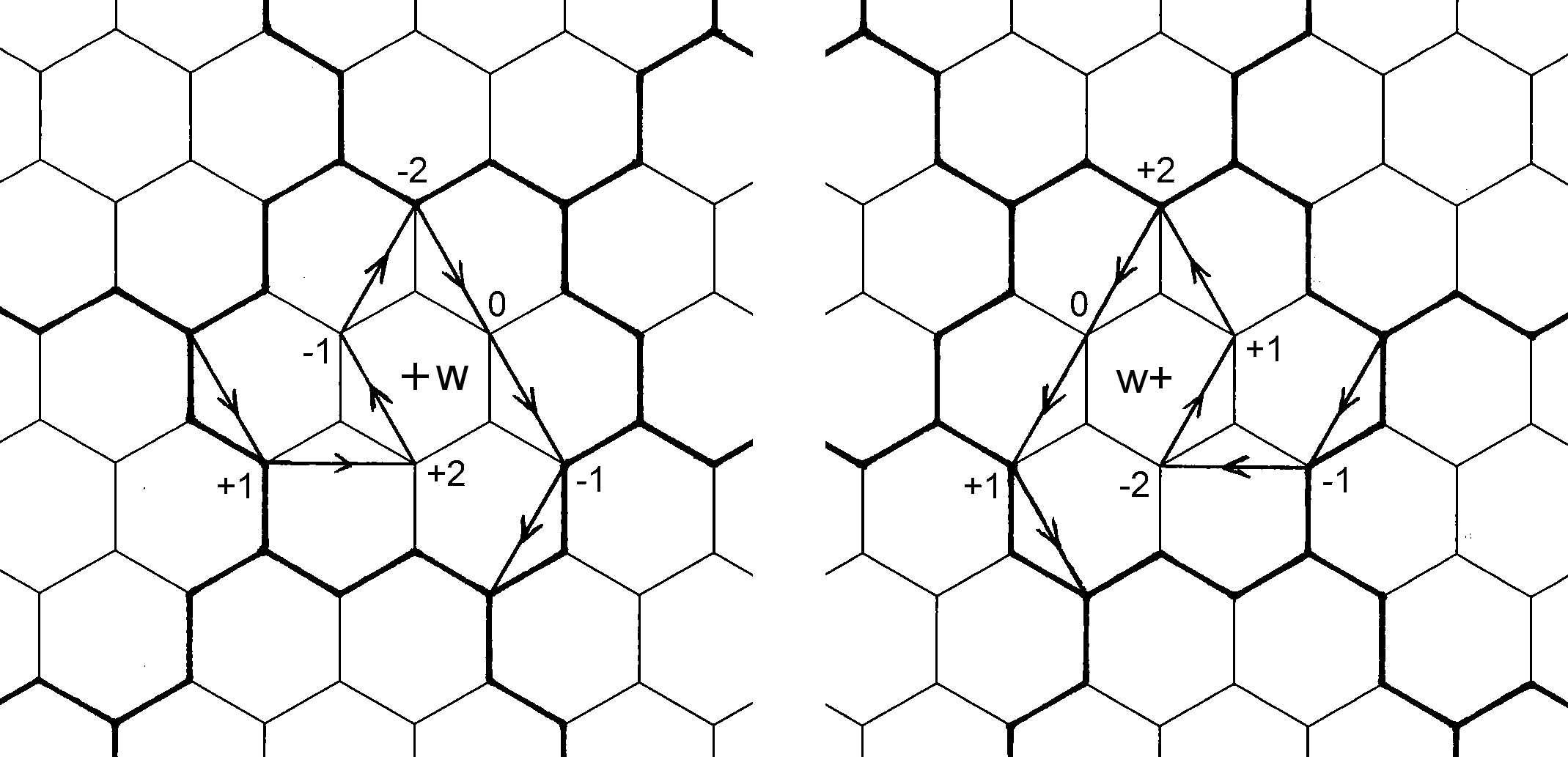}

\medskip Figure 3.2
\end{center}

The first (resp. second) part of Figure 3.3 below shows a covering of a tile
of $\mathcal{P}_{w++}$ (resp. $\mathcal{P}_{w+-}$) by a curve associated to $%
S_{++}$ (resp. $S_{+-}$).

\begin{center}
\includegraphics[scale=0.72]{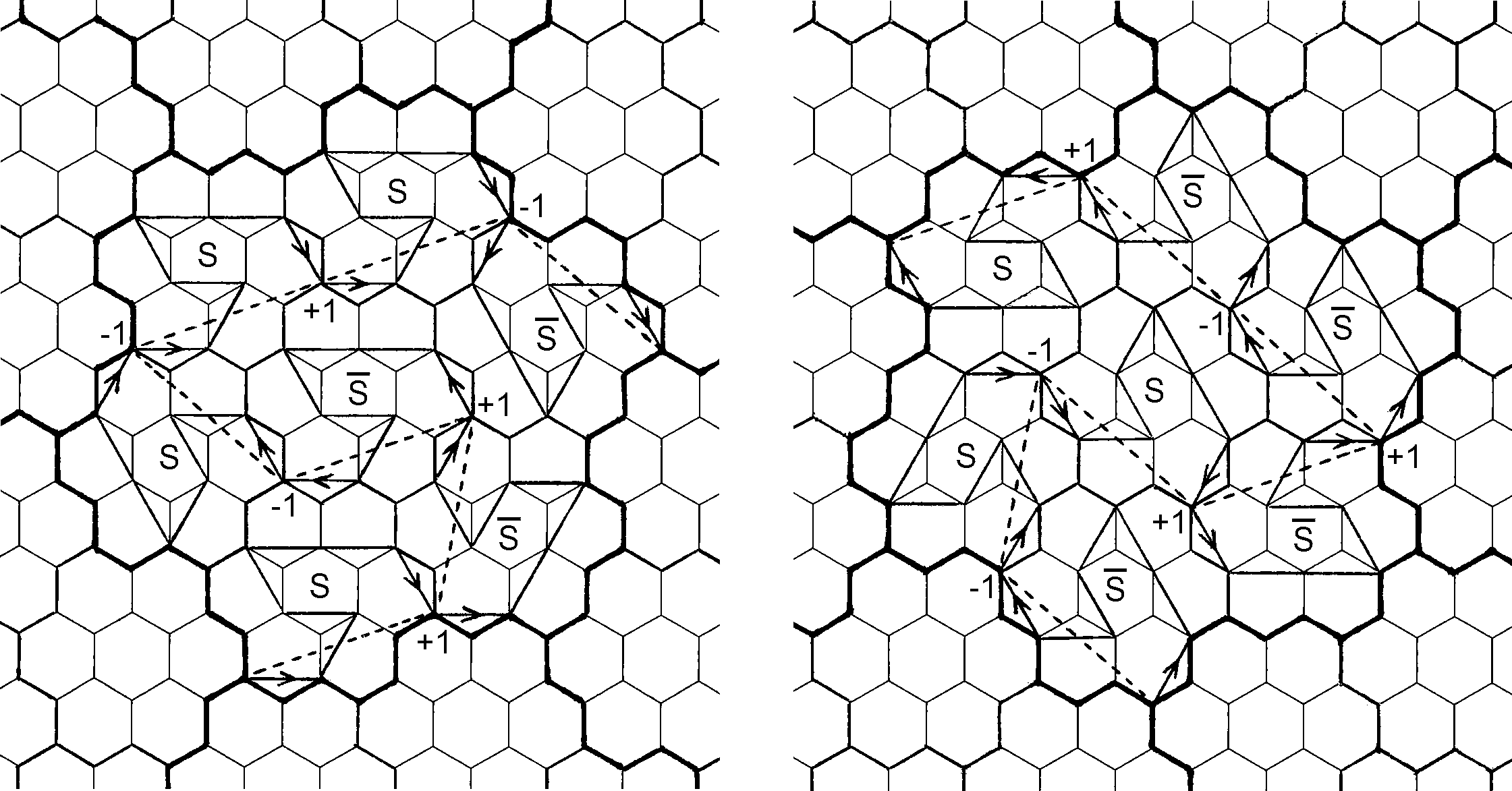}

\medskip Figure 3.3
\end{center}

\noindent \textbf{Proposition 3.2.} Consider $w\in 
%TCIMACRO{\U{211d} }%
%BeginExpansion
\mathbb{R}
%EndExpansion
^{2}$, $n\in 
%TCIMACRO{\U{2115} }%
%BeginExpansion
\mathbb{N}
%EndExpansion
^{\ast }$,\ $\lambda _{1},\ldots ,\lambda _{n}\in \left\{ -,+\right\} $ and $%
P\in \mathcal{P}_{w\lambda _{1}\cdots \lambda _{n}}$. Then $P$ has $6$\
coverings by oriented curves, each of them determined by its initial and
terminal points. Each covering of $P$ is associated to $S_{\lambda
_{1}\cdots \lambda _{n}}$ or $\overline{S_{\lambda _{1}\cdots \lambda _{n}}}$%
. The $3$ coverings associated to $S_{\lambda _{1}\cdots \lambda _{n}}$ are
obtained from one of them by rotations of angles $2k\pi /3$, and the $3$
others by changing their orientations. For $n\geq 2$, if $C$ is a covering
of $P$ associated to $S_{\lambda _{1}\cdots \lambda _{n}}$ (resp. $\overline{%
S_{\lambda _{1}\cdots \lambda _{n}}}$), then $\Delta _{w\lambda _{1}}(C)$ is
a covering of $\Delta _{w\lambda _{1}}(C)$ associated to $S_{\lambda
_{2}\cdots \lambda _{n}}$ (resp. $\overline{S_{\lambda _{2}\cdots \lambda
_{n}}}$).\bigskip

\noindent \textbf{Proof.} We see from Figure 3.2 that Proposition 3.2 is
true for $n=1$. Now we show that it is true for $n\geq 2$ if it is true for $%
n-1$.

By the induction hypothesis, there exists a covering $D$ of $\Delta
_{w\lambda _{1}}(P)$ associated to $S_{\lambda _{2}\cdots \lambda _{n}}$.
Each segment $A$\ of $D$\ is a covering of an hexagon $H_{A}\in \Delta
_{w\lambda _{1}}(\mathcal{P}_{w\lambda _{1}\cdots \lambda _{n}})$\ which is
contained in $\Delta _{w\lambda _{1}}(P)$. We obtain a covering $C$ of $P$
with $\Delta _{w\lambda _{1}}(C)=D$ by replacing each such $A$\ with a
covering of $\Delta _{w\lambda _{1}}^{-1}(H_{A})$ which has the
corresponding initial and terminal points.

As $P$ is invariant through a rotation of angle $2\pi /3$, it has $3$
coverings obtained from $C$ by rotations of angles $0,2\pi /3,4\pi /3$, and $%
3$ others obtained by changing their orientations. These coverings are
associated to the $6$ possible pairs of initial and terminal points.

Any other\ covering $B$ of $P$ has the same initial and terminal points as
one of the $6$ coverings above, say $C^{\prime }$. Then $\Delta _{w\lambda
_{1}}(B)$ and $\Delta _{w\lambda _{1}}(C^{\prime })$ are coverings of $%
\Delta _{w\lambda _{1}}(P)$ with the same initial and terminal points. By
the induction hypothesis, we have $\Delta _{w\lambda _{1}}(B)=\Delta
_{w\lambda _{1}}(C^{\prime })$. Consequently, for each $Q\in \mathcal{P}%
_{w\lambda _{1}}$ contained in $P$, $B\upharpoonright Q$\ and $C^{\prime
}\upharpoonright Q$\ are coverings of $Q$ with the same initial and terminal
points, which implies $B\upharpoonright Q=C^{\prime }\upharpoonright Q$. It
follows $B=C^{\prime }$.

Now it suffices to show that $C$ is associated to $S_{\lambda _{1}\cdots
\lambda _{n}}$. We see from Figure 3.3 that it is true for $n=2$. For $n\geq
3$, we write $C=(C_{h})_{0\leq h\leq 7^{n}-1}$ and $D=(D_{h})_{0\leq h\leq
7^{n-1}-1}$. We consider the associated sequences

\noindent $(c_{h})_{1\leq h\leq 7^{n}-1}\in \left\{ -2,-1,0,+1,+2\right\}
^{7^{n}-1}$ and

\noindent $(d_{h})_{1\leq h\leq 7^{n-1}-1}\in \left\{ -2,-1,0,+1,+2\right\}
^{7^{n-1}-1}$.

For $0\leq k\leq 7^{n-2}-1$, we write $U_{k}=(c_{h})_{49k+1\leq h\leq
49k+48} $\ and $V_{k}=(d_{h})_{7k+1\leq h\leq 7k+6}$. There exists $Q_{k}\in 
\mathcal{P}_{w\lambda _{1}\lambda _{2}}$\ contained in $P$ such that $U_{k}$
is associated to $C\upharpoonright Q_{k}$ and $V_{k}$ is associated to $%
D\upharpoonright \Delta _{w\lambda _{1}}(Q_{k})=\Delta _{w\lambda
_{1}}(C\upharpoonright Q_{k})$. Consequently, we have $U_{k}=S_{\lambda
_{1}\lambda _{2}}$ (resp. $\overline{S_{\lambda _{1}\lambda _{2}}}$) if and
only if $V_{k}=S_{\lambda _{2}}$ (resp. $\overline{S_{\lambda _{2}}}$).

It remains to be proved that $c_{49k}=d_{7k}$ for $1\leq k\leq 7^{n-2}-1$.\
As $U_{k-1},U_{k}\in \{S_{\lambda _{1}\lambda _{2}},\overline{S_{\lambda
_{1}\lambda _{2}}}\}$, there exists $T\in \{S_{\lambda _{1}},\overline{%
S_{\lambda _{1}}}\}$ such that $U_{k-1}$ is ending with $T$ and $U_{k}$ is
beginning with $\overline{T}$. It follows that the angle between $C_{49k}$
and $C_{49k+1}$ is equal to the angle between $D_{7k}$ and $D_{7k+1}$, which
implies $c_{49k}=d_{7k}$.~~$\blacksquare $\bigskip

\noindent \textbf{Example}. The Peano-Gosper curves considered in [3, p. 46]
and [2, p. 63] are associated to the sequences $T_{n}=S_{\lambda _{1}\cdots
\lambda _{n}}$ with $\lambda _{1}=\cdots =\lambda _{n}=+$. The case $n=2$ is
shown in the first part of Figure 3.3. We have

\noindent $T_{n+1}=(T_{n},+1,\overline{T_{n}},+1,\overline{T_{n}}%
,-1,T_{n},-1,T_{n},+1,T_{n},-1,\overline{T_{n}})$

\noindent for each $n\in 
%TCIMACRO{\U{2115} }%
%BeginExpansion
\mathbb{N}
%EndExpansion
^{\ast }$. In [2] and several papers mentioned among its references, W.
Kuhirun, D.H. Werner and P.L. Werner prove that an antenna with the shape of
a Peano-Gosper curve has particular electromagnetic properties. We can
suppose that similar properties exist for the other values of $\lambda
_{1},\ldots ,\lambda _{n}$.\bigskip

\noindent \textbf{Corollary 3.3.} For each $w\in 
%TCIMACRO{\U{211d} }%
%BeginExpansion
\mathbb{R}
%EndExpansion
^{2}$, each $n\in 
%TCIMACRO{\U{2115} }%
%BeginExpansion
\mathbb{N}
%EndExpansion
$, any\ $\lambda _{1},\ldots ,\lambda _{n+1}\in \left\{ -,+\right\} $
(resp.\ $\lambda _{1},\ldots ,\lambda _{n+2}\in \left\{ -,+\right\} $) and
each $P\in \mathcal{P}_{w\lambda _{1}\cdots \lambda _{n}}$, each covering of
some $Q\in \mathcal{P}_{w\lambda _{1}\cdots \lambda _{n+1}}$ (resp. $%
\mathcal{P}_{w\lambda _{1}\cdots \lambda _{n+2}}$) by a nonoriented (resp.
oriented) curve contains copies of the $3$ (resp. $6$) coverings of $P$ by
nonoriented (resp. oriented) curves.\bigskip

\noindent \textbf{Proof}. By Proposition 3.2, it suffices to show that each
covering of some $Q\in \mathcal{P}_{w\lambda _{1}\cdots \lambda _{n+1}}$
(resp. $\mathcal{P}_{w\lambda _{1}\cdots \lambda _{n+2}}$) by a nonoriented
(resp. oriented) curve contains $3$ (resp. $6$) nonisomorphic coverings of
tiles $P\in \mathcal{P}_{w\lambda _{1}\cdots \lambda _{n}}$ by nonoriented
(resp. oriented) curves. We see from Figure 3.2 (resp. 3.3) that the
statement for nonoriented (resp. oriented) curves is true for $n=0$.

For each $n\geq 1$ and each covering $C$ of some $Q\in \mathcal{P}_{w\lambda
_{1}\cdots \lambda _{n+1}}$ (resp. $\mathcal{P}_{w\lambda _{1}\cdots \lambda
_{n+2}}$) by a nonoriented (resp. oriented) curve, we consider some
hexagonal tiles $Q_{1},Q_{2},Q_{3}$ (resp. $Q_{1},\ldots ,Q_{6}$) contained
in $\Delta _{w\lambda _{1}\cdots \lambda _{n}}(Q)$ such that the nonoriented
(resp. oriented) segments $\Delta _{w\lambda _{1}\cdots \lambda
_{n}}(C)\upharpoonright Q_{i}$ are nonisomorphic. Then the nonoriented
(resp. oriented) curves $C\upharpoonright \Delta _{w\lambda _{1}\cdots
\lambda _{n}}^{-1}(Q_{i})=\Delta _{w\lambda _{1}\cdots \lambda
_{n}}^{-1}(\Delta _{w\lambda _{1}\cdots \lambda _{n}}(C)\upharpoonright
Q_{i})$\ are nonisomorphic.~~$\blacksquare $\bigskip

\noindent \textbf{Lemma 3.4.} Here the curves are nonoriented. Consider $%
w\in 
%TCIMACRO{\U{211d} }%
%BeginExpansion
\mathbb{R}
%EndExpansion
^{2}$, $n\in 
%TCIMACRO{\U{2115} }%
%BeginExpansion
\mathbb{N}
%EndExpansion
$,\ $\lambda _{1},\ldots ,\lambda _{n+1}\in \left\{ -,+\right\} $, $P\in 
\mathcal{P}_{w\lambda _{1}\cdots \lambda _{n+1}}$,\ $Q_{1},\ldots ,Q_{7}\in 
\mathcal{P}_{w\lambda _{1}\cdots \lambda _{n}}$ such that $P=Q_{1}\cup
\cdots \cup Q_{7}$, and $Q\in \left\{ Q_{1},\ldots ,Q_{7}\right\} $. If $Q$
contains the center of $P$, then each of the $3$ coverings of $Q$ extends
into a covering of $P$. Otherwise, denote by $x_{1}$ (resp. $x_{2}$, $x_{3}$%
) the vertex of $Q$\ belonging to $W$ which is a vertex of $1$ (resp. $2$, $%
3 $) tiles\ $Q_{i}$.\ If $x_{1}$ is not a vertex of $P$, then the covering
of $Q$ which joins $x_{2}$ and $x_{3}$ extends into $3$ coverings of $P$ and
the $2$ other coverings of $Q$ do not extend.\ If $x_{1}$ is a vertex of $P$%
, then the covering of $Q$ which joins $x_{1}$ and $x_{2}$ (resp. $x_{2}$
and $x_{3}$, $x_{1}$ and $x_{3}$) extends into $2$ (resp. $1$, $0$)
coverings of $P$.\bigskip

\noindent \textbf{Proof}. We see from Figure 3.2 that Lemma 3.4 is true for $%
n=0$. Then, as in the proof of Corollary 3.3, we use the derivations $\Delta
_{w\lambda _{1}\cdots \lambda _{n}}$\ to show that it is also true for $%
n\geq 1$.~~$\blacksquare $\bigskip

For each $\mathcal{P}_{X\Lambda }$, we say that a curve $C$ is a \emph{%
covering} of a region $R$ of $\mathcal{P}_{X\Lambda }$ if $C\subset R$ and
if $C$ contains a covering of each $P\in \mathcal{P}_{X\Lambda }$ contained
in $R$. We say that a set $\mathcal{C}$\ of disjoint complete curves is a 
\emph{covering} of $\mathcal{P}_{X\Lambda }$ if each $P\in \mathcal{P}%
_{X\Lambda }$\ has a covering by a subcurve of a curve of $\mathcal{C}$.
Then each region of $\mathcal{P}_{X\Lambda }$\ also has a covering by a
subcurve of a curve of $\mathcal{C}$.\bigskip

\noindent \textbf{Proposition 3.5.} We consider nonoriented curves. Any
region $R$ of some $\mathcal{P}_{X\Lambda }$ has a covering by a complete
curve. This covering is unique if $R\varsubsetneq 
%TCIMACRO{\U{211d} }%
%BeginExpansion
\mathbb{R}
%EndExpansion
^{2}$. If $R$ is the union of an increasing sequence of tiles with $1$
common vertex $y\in W$, then $R$ also has a unique covering by a half curve\
and $y$ is the endpoint of that curve. Otherwise, no such covering
exists.\bigskip

\noindent \textbf{Proof}. We have $R=\cup _{n\in 
%TCIMACRO{\U{2115} }%
%BeginExpansion
\mathbb{N}
%EndExpansion
}P_{n}$ for an increasing sequence $(P_{n})_{n\in 
%TCIMACRO{\U{2115} }%
%BeginExpansion
\mathbb{N}
%EndExpansion
}\in \times _{n\in 
%TCIMACRO{\U{2115} }%
%BeginExpansion
\mathbb{N}
%EndExpansion
}\mathcal{P}_{X\Lambda }^{n}$.

For each $m\in 
%TCIMACRO{\U{2115} }%
%BeginExpansion
\mathbb{N}
%EndExpansion
$, consider the nonempty set $E_{m}$ which consists of the coverings of $%
P_{m}$ with endpoints different from the common vertex of the tiles $P_{n}$
if it exists. Then, for $m<n$, each element of $E_{n}$\ gives by restriction
an element of $E_{m}$. By K\"{o}nig's lemma, there exists an increasing
sequence $(C_{n})_{n\in 
%TCIMACRO{\U{2115} }%
%BeginExpansion
\mathbb{N}
%EndExpansion
}\in \times _{n\in 
%TCIMACRO{\U{2115} }%
%BeginExpansion
\mathbb{N}
%EndExpansion
}E_{n}$; $\cup _{n\in 
%TCIMACRO{\U{2115} }%
%BeginExpansion
\mathbb{N}
%EndExpansion
}C_{n}$\ is a covering of $R$ by a complete curve.

Now we can suppose $R\varsubsetneq 
%TCIMACRO{\U{211d} }%
%BeginExpansion
\mathbb{R}
%EndExpansion
^{2}$ for the remainder of the proof. Then there exists $h\in 
%TCIMACRO{\U{2115} }%
%BeginExpansion
\mathbb{N}
%EndExpansion
$ such that, for each $n\geq h$, $P_{n}$ does not contain the center of $%
P_{n+1}$. For $n\geq h$, denote by $y_{n}^{1}$ (resp. $y_{n}^{2}$,~$%
y_{n}^{3} $)\ the vertex of $P_{n}$\ belonging to $W$ which is a vertex of $%
1 $ (resp. $2$,~$3$) tiles of $\mathcal{P}_{X\Lambda }^{n}$ contained in $%
P_{n+1}$.

If $R$ is the union of an increasing sequence of tiles with a common vertex
belonging to $W$, then there exists $k\geq h$ such that $y_{n}^{1}=y_{k}^{1}$
for each $n\geq k$. For $n\geq k$, denote by $F_{n}$ the nonempty set which
consists of the coverings of $P_{n}$\ with endpoint $y_{k}^{1}$.\ Then, for $%
k\leq m\leq n$, each element of $F_{n}$ gives by restriction an element of $%
F_{m}$. By K\"{o}nig's lemma, there exists an increasing sequence $%
(D_{n})_{n\geq k}\in \times _{n\geq k}F_{n}$; $\cup _{n\in 
%TCIMACRO{\U{2115} }%
%BeginExpansion
\mathbb{N}
%EndExpansion
}D_{n}$\ is a covering of $R$ by a half curve with endpoint $y_{k}^{1}$.

Conversely, suppose that $R$ has $2$ coverings $C,D$. For each $n\in 
%TCIMACRO{\U{2115} }%
%BeginExpansion
\mathbb{N}
%EndExpansion
$, write $C_{n}=C\upharpoonright P_{n}$ and $D_{n}=D\upharpoonright P_{n}$.
Consider $k\geq h$ such that $C_{k}\neq D_{k}$.

For $n\geq k$, we have $C_{n}\neq D_{n}$. As $C_{n}=C_{n+1}\upharpoonright
P_{n}$ and $D_{n}=D_{n+1}\upharpoonright P_{n}$, it follows from Lemma 3.4
applied to $P_{n}$ and $P_{n+1}$ that $y_{n}^{1}$ is a vertex of $P_{n+1}$
and that one of the curves $C_{n},D_{n}$ connects $y_{n}^{1}$ and $y_{n}^{2}$%
, while the other one connects $y_{n}^{2}$ and $y_{n}^{3}$.

Consequently, we have $y_{n}^{1}=y_{k}^{1}$ for $n\geq k$, one of the curves 
$C,D$ is a complete curve, the other one is a half curve with endpoint $%
y_{k}^{1}$ and there is only one possibility for each of them.~~$%
\blacksquare $\bigskip

Now we state and prove our main results. As above, we only consider curves
with segments in $W$.\bigskip

\noindent \textbf{Theorem 3.6.} For each $\Lambda $, we consider the
coverings of the sets $\mathcal{P}_{X\Lambda }$ by sets of nonoriented
complete curves:

\noindent 1) Suppose that $\mathcal{P}_{X\Lambda }$\ has $1$ region. Then
each covering of $\mathcal{P}_{X\Lambda }$ consists of $1$ curve. If $X$ is
ultimately constant, then $\mathcal{P}_{X\Lambda }$ has $3$ coverings,
obtained from one another by rotations of angles $\pm 2\pi /3$. Otherwise, $%
\mathcal{P}_{X\Lambda }$ has $1$ or $2$ coverings, and each case is realized
for $2^{\omega }$ values of $X$.

\noindent 2) If $\mathcal{P}_{X\Lambda }$ has $2$ regions, then $\mathcal{P}%
_{X\Lambda }$ has $1$ covering. It consists of $2$ curves. This case is
realized for $2^{\omega }$ values of $X$.

\noindent 3) If $\mathcal{P}_{X\Lambda }$ has $3$ regions and if their
common point $y$ does not belong to $W$, then $\mathcal{P}_{X\Lambda }$ has $%
1$ covering. It consists of $3$ curves obtained from one another by
rotations of center $y$ and angles $\pm 2\pi /3$.

\noindent 4) If $\mathcal{P}_{X\Lambda }$ has $3$ regions and if their
common point $y$ belongs to $W$, then $\mathcal{P}_{X\Lambda }$ has $3$
coverings obtained from one another by rotations of center $y$ and angles $%
\pm 2\pi /3$. Each covering of $\mathcal{P}_{X\Lambda }$ consists of $2$
curves. One of them is a covering of a region. The other one is the union of 
$2$ half curves with endpoint $y$, which are equivalent modulo a rotation of
center $y$ and angle $2\pi /3$; each half curve is a covering of a
region.\bigskip

Theorem 3.6 will be proved after the remarks and the example below, which
concern nonoriented curves:\bigskip 

\noindent \textbf{Remark.} It follows from Theorem 3.6 that, for each region 
$R\varsubsetneq 
%TCIMACRO{\U{211d} }%
%BeginExpansion
\mathbb{R}
%EndExpansion
^{2}$ of some $\mathcal{P}_{X\Lambda }$, the covering of $R$ by a complete
curve can be extended into $1$ covering of $\mathcal{P}_{X\Lambda }$. The
covering of $R$ by a half curve, if it exists, can be extended into $2$
coverings of $\mathcal{P}_{X\Lambda }$ by complete curves, which are
equivalent modulo a rotation of angle $2\pi /3$.\bigskip

\noindent \textbf{Remark}. For each $\Lambda $, as each $\mathcal{P}%
_{X\Lambda }$ has finitely many coverings, Proposition 3.1 implies that each
isomorphism class of coverings of sets $\mathcal{P}_{X\Lambda }$ is
countable. Consequently, it follows from Theorem 3.6 that we have $2^{\omega
}$ isomorphism classes of coverings for case 2) and $2^{\omega }$ classes
for case 1) with $X$ not ultimately constant. On the other hand, Proposition
3.1 and Theorem 3.6 imply that we have $1$ isomorphism class for case 3), $3$
classes for case 4) and $3$ classes for case 1) with $X$ ultimately constant.

\begin{center}
\includegraphics[scale=0.66]{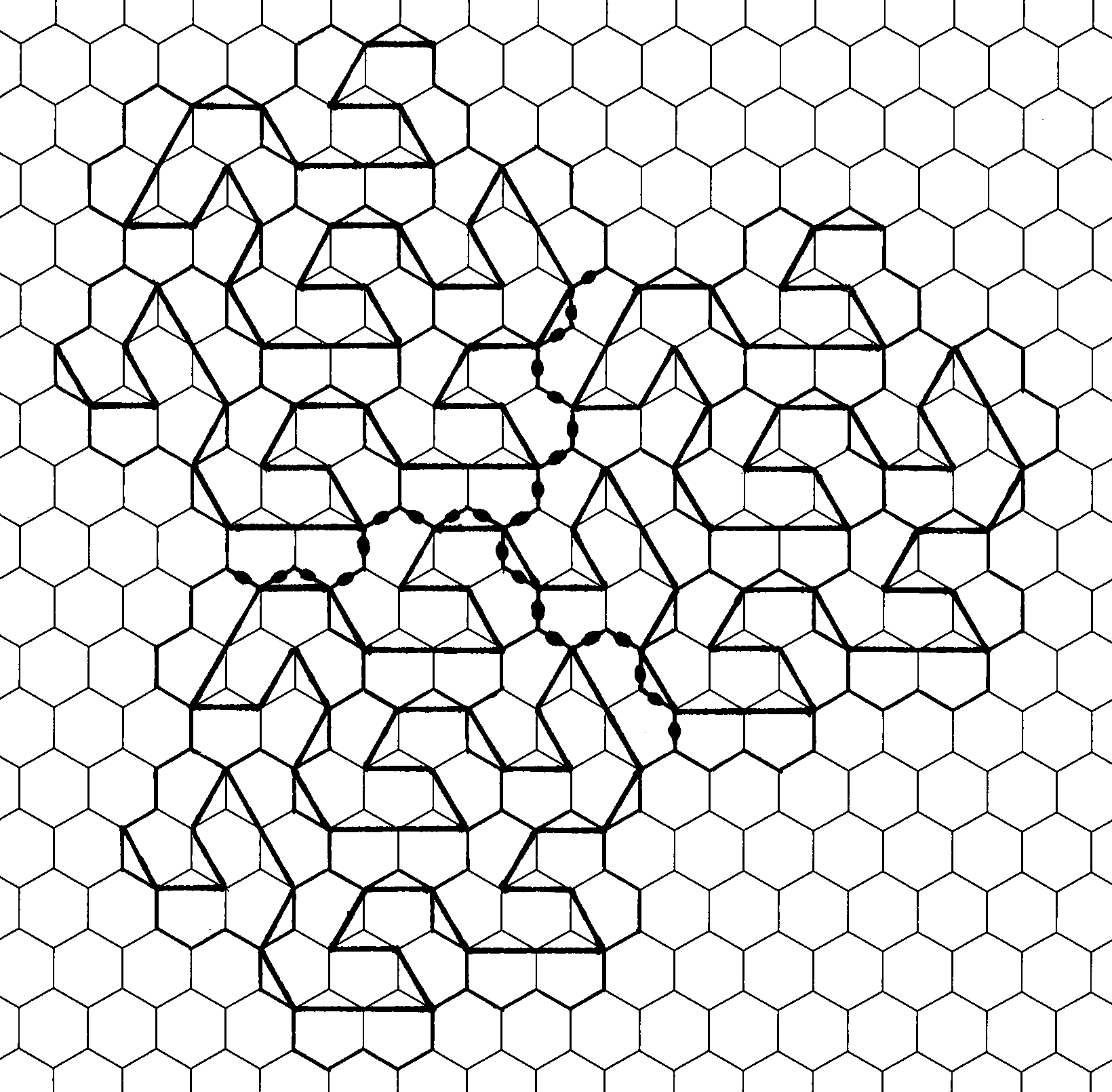}\medskip

Figure 3.4
\end{center}

\noindent \textbf{Example}. For each $n\in 
%TCIMACRO{\U{2115} }%
%BeginExpansion
\mathbb{N}
%EndExpansion
^{\ast }$, denote by $C_{n}$ the Peano-Gosper curve associated to the
sequence $T_{n}$. Consider the inductive limit $B$ (resp. $C,D$), defined up
to rotation, of the curves $C_{n}$ with each $C_{n}$ embedded in $C_{n+1}$
by identifying $T_{n}$ with its fourth (resp. third, first) copy in

\noindent $T_{n+1}=(T_{n},+1,\overline{T_{n}},+1,\overline{T_{n}}%
,-1,T_{n},-1,T_{n},+1,T_{n},-1,\overline{T_{n}})$.

\noindent Then there exist: 1) a covering of $%
%TCIMACRO{\U{211d} }%
%BeginExpansion
\mathbb{R}
%EndExpansion
^{2}$ by $3$ copies of the complete curve $B$; 2) a covering of $%
%TCIMACRO{\U{211d} }%
%BeginExpansion
\mathbb{R}
%EndExpansion
^{2}$ by a copy of the complete curve $C$, and $2$ copies of the half curve $%
D$ which form a complete curve (see Figure 3.4 above).\bigskip

\noindent \textbf{Proof of Theorem 3.6.} We write $X=(x_{n})_{n\in 
%TCIMACRO{\U{2115} }%
%BeginExpansion
\mathbb{N}
%EndExpansion
}$. For each$\ n\in 
%TCIMACRO{\U{2115} }%
%BeginExpansion
\mathbb{N}
%EndExpansion
$, we consider the tile $P_{n}\in \mathcal{P}_{X\Lambda }^{n}$ which
contains $x_{n}$.

\noindent \textbf{Proof of 1)}. We have $%
%TCIMACRO{\U{211d} }%
%BeginExpansion
\mathbb{R}
%EndExpansion
^{2}=\cup _{n\in 
%TCIMACRO{\U{2115} }%
%BeginExpansion
\mathbb{N}
%EndExpansion
}P_{n}$. Each covering of $\mathcal{P}_{X\Lambda }$ consists of $1$ complete
curve by Proposition 3.5.

First we suppose that there exists $h$ such that $x_{n}=x_{h}$ for each $%
n\geq h$.\ Then, by Proposition 3.2, for $m\geq h$,\ the $3$ coverings of $%
P_{m}$ are\ obtained from one of them by rotations of center $x_{h}$ and
angles $2k\pi /3$, and they are the\ restrictions of\ the $3$ coverings of $%
P_{n}$ for each $n\geq m$. It follows that $\mathcal{P}_{X\Lambda }$ has
exactly $3$ coverings, which are the inductive limits of\ the coverings of
the tiles $P_{n}$ for $n\geq h$, and therefore\ obtained from one of them by
rotations of center $x_{h}$ and angles $2k\pi /3$.

From now on, we suppose $X$ not ultimately constant. First we observe that $%
\mathcal{P}_{X\Lambda }$ cannot have $3$ coverings. Otherwise, for $n$ large
enough, these $3$ coverings would give by restriction $3$ distinct coverings
of $P_{n}$, which is not possible when $x_{n+1}\neq x_{n}$ since, by Lemma
3.4, at most $2$ coverings of $P_{n}$ extend into coverings of $P_{n+1}$.

By Lemma 3.4, for each $n\in 
%TCIMACRO{\U{2115} }%
%BeginExpansion
\mathbb{N}
%EndExpansion
$ and any coverings $A\neq B$ of $P_{n}$, we have $2$ coverings of $P_{n+1}$
which extend $A$ and $B$ in $2$ cases: first if $P_{n}$ contains the center
of $P_{n+1}$, second if $P_{n}$ does not contain the center of $P_{n+1}$, if
the common endpoint of $A$ and $B$ belongs to exactly $2$ tiles of $\mathcal{%
P}_{X\Lambda }^{n}$ contained in $P_{n+1}$, and if the other endpoint of $A$
or the other endpoint of $B$ is a vertex of $P_{n+1}$. The first case is
realized for $x_{n+1}=x_{n}$ and the second one for another choice of $%
x_{n+1}$.

It follows that, for any coverings $C\neq D$ of some $P\in \mathcal{P}$ and
each $K\subset 
%TCIMACRO{\U{2115} }%
%BeginExpansion
\mathbb{N}
%EndExpansion
$ with $K$ and $%
%TCIMACRO{\U{2115} }%
%BeginExpansion
\mathbb{N}
%EndExpansion
-K$ infinite, we can choose $X=(x_{n})_{n\in 
%TCIMACRO{\U{2115} }%
%BeginExpansion
\mathbb{N}
%EndExpansion
}$, with $x_{0}$ center of $P$\ and $x_{n+1}=x_{n}$ if and only if $n\in K$,
so that there exist some increasing sequences $(C_{n})_{n\in 
%TCIMACRO{\U{2115} }%
%BeginExpansion
\mathbb{N}
%EndExpansion
}$, $(D_{n})_{n\in 
%TCIMACRO{\U{2115} }%
%BeginExpansion
\mathbb{N}
%EndExpansion
}$ with $C_{0}=C$, $D_{0}=D$ and $C_{n},D_{n}$ coverings of $P_{n}$ for each 
$n\in 
%TCIMACRO{\U{2115} }%
%BeginExpansion
\mathbb{N}
%EndExpansion
$. For each $n\in 
%TCIMACRO{\U{2115} }%
%BeginExpansion
\mathbb{N}
%EndExpansion
-K$, supposing that $x_{n},C_{n},D_{n}$ are already defined, we choose $%
x_{n+1}$ in such a way that the second case above is realized for $A=C_{n}$
and $B=D_{n}$. Then $\mathcal{P}_{X\Lambda }$ only has $1$ region and $\cup
_{n\in 
%TCIMACRO{\U{2115} }%
%BeginExpansion
\mathbb{N}
%EndExpansion
}C_{n}$, $\cup _{n\in 
%TCIMACRO{\U{2115} }%
%BeginExpansion
\mathbb{N}
%EndExpansion
}D_{n}$ are the only coverings of $\mathcal{P}_{X\Lambda }$.

By Lemma 3.4, for each $n\in 
%TCIMACRO{\U{2115} }%
%BeginExpansion
\mathbb{N}
%EndExpansion
$, if $x_{n+1}\neq x_{n}$ and if the point of $P_{n}$ which is a vertex of $%
P_{n+1}$ does not belong to $W$, then the $3$ coverings of $P_{n+1}$ give by
restriction the same covering of $P_{n}$, so that only $1$ covering of $%
P_{n} $ is the restriction of a covering of $\mathcal{P}_{X\Lambda }$. There
are $2^{\omega }$ different ways to choose $X$ so that this property is true
for infinitely many integers $n$, and so that $x_{n+1}=x_{n}$ is also true
for infinitely many integers $n$. Then $\mathcal{P}_{X\Lambda }$ only has $1$
covering and $1$ region.

\noindent \textbf{Proof of 2), 3), 4)}. Suppose that $\mathcal{P}_{X\Lambda
} $ has $2$ or $3$\ regions. Write $R=\cup _{n\in 
%TCIMACRO{\U{2115} }%
%BeginExpansion
\mathbb{N}
%EndExpansion
}P_{n}$.

If $R$ is not the union of an increasing sequence of tiles with $1$ common
vertex, then $\mathcal{P}_{X\Lambda }$ has $2$\ regions and the second
region satisfies the same property. By Proposition 3.5, each region has a
unique covering and it is a complete curve. The union of these coverings is
the unique covering of $\mathcal{P}_{X\Lambda }$.

This case is realized if:

\noindent a) there exists a sequence $(\Sigma _{n})_{n\in 
%TCIMACRO{\U{2115} }%
%BeginExpansion
\mathbb{N}
%EndExpansion
}$\ with $\Sigma _{n}$\ side of $P_{n}$\ and $\Sigma _{n}\subset \Sigma
_{n+1}$\ for each\ $n\in 
%TCIMACRO{\U{2115} }%
%BeginExpansion
\mathbb{N}
%EndExpansion
$;

\noindent b) there is no common vertex of the tiles $P_{n}$\ for $n$ large.

For each\ $n\in 
%TCIMACRO{\U{2115} }%
%BeginExpansion
\mathbb{N}
%EndExpansion
$\ and each choice of $x_{1},\ldots ,x_{n}$ compatible with a), there are $3$
choices of $x_{n+1}$ compatible with a). Consequently, there are $2^{\omega
} $ sequences $X$ which satisfy a), and also $2^{\omega }$ sequences which
satisfy a) and b) since countably many sequences do not satisfy b).

If $R$ is the union of an increasing sequence of tiles with one common
vertex $y$, then $\mathcal{P}_{X\Lambda }$ has $3$\ regions $%
R_{1},R_{2},R_{3}$ obtained from $R$\ by rotations of center $y$ and angles $%
2k\pi /3$.

If $y$\ does not belong to $W$, then, by Proposition 3.5, each $R_{i}$ has a
unique covering. These coverings are complete curves obtained from one of
them\ by rotations of center $y$ and angles $2k\pi /3$. Their union is the
unique covering of $\mathcal{P}_{X\Lambda }$.

If $y$ belongs to $W$, then, by Proposition 3.5, each $R_{i}$ has $1$
covering by a complete curve and $1$ covering by a half curve with endpoint $%
y$. Each of the $3$\ coverings of $\mathcal{P}_{X\Lambda }$ is obtained by
taking a covering by a complete curve for $1$ region and $2$ coverings by
half curves for the $2$ other regions. The $2$ half curves are equivalent
modulo a rotation of center $y$ and angle $2\pi /3$. They form a complete
curve since $y$ is their common endpoint.~~$\blacksquare $\bigskip

Now, for each set $\mathcal{C}$ of oriented curves, we consider the
following property:

\noindent (P) If $2$ segments of curves of $\mathcal{C}$ are opposite sides
of a rhombus, then they have opposite orientations.

\noindent We observe that (P) is satisfied if $\mathcal{C}$ consists of $1$
complete curve or $1$ curve which is a covering of a tile of some $\mathcal{P%
}_{X\Lambda }$.\bigskip

\noindent \textbf{Theorem 3.7.} The local isomorphism property implies (P)
for coverings by sets of oriented curves. The strong local isomorphism
property is true for the coverings of sets $\mathcal{P}_{X\Lambda }$ by sets
of nonoriented curves or by sets of oriented curves which satisfy
(P).\bigskip

\noindent \textbf{Remark}. Theorem 3.7 implies that each covering of some $%
\mathcal{P}_{X\Lambda }$ by a set of oriented curves satisfies the local
isomorphism property or can be transformed into a covering which satisfies
that property by changing the orientation of one of the curves.\bigskip

\noindent \textbf{Theorem 3.8.} For each $\Lambda $\ and any $X,Y$, any\
coverings of $\mathcal{P}_{X\Lambda }$ and $\mathcal{P}_{Y\Lambda }$ by sets
of nonoriented curves are locally isomorphic. Any\ coverings of $\mathcal{P}%
_{X\Lambda }$ and $\mathcal{P}_{Y\Lambda }$ by sets of oriented curves are
locally isomorphic if they satisfy (P).\bigskip

\noindent \textbf{Remark}. It follows from Theorem 3.6 that each local
isomorphism class of coverings by sets of nonoriented curves associated to a
sequence $\Lambda $ contains $2^{\omega }$\ (resp. $2^{\omega }$, $1$, $0$)\
isomorphism classes of coverings by $1$ (resp. $2$, $3$, $\geq 4$)
curves.\bigskip

\noindent \textbf{Proof of Theorem 3.7 and Theorem 3.8}. It suffices to
prove the results for oriented curves since they imply the results for
nonoriented curves.

First suppose that a covering $\mathcal{C}$ of some $\mathcal{P}_{X\Lambda }$
satisfies the local isomorphism property. Consider $2$ segments of curves of 
$\mathcal{C}$ which are opposite sides of a rhombus $A$. Then there exist a
translation $\tau $\ and a tile $T$ of $\mathcal{P}_{X\Lambda }$ such that $%
\tau (A)$\ is contained in $T$ and $\tau (A\cap \mathcal{C})=\tau (A)\cap C$%
, where $C$ is the curve of $\mathcal{C}$ which contains a coverig of $T$.
The $2$ segments have opposite orientations because their images through $%
\tau $, which belong to the same curve, necessarily have opposite
orientations.

Now, for each $X$ such that $\mathcal{P}_{X\Lambda }$ exists, each $n\in 
%TCIMACRO{\U{2115} }%
%BeginExpansion
\mathbb{N}
%EndExpansion
$ and each $x\in 
%TCIMACRO{\U{211d} }%
%BeginExpansion
\mathbb{R}
%EndExpansion
^{2}$, we consider the tile $P_{Xx}^{n}\in \mathcal{P}_{X\Lambda }^{n}$\
with center $x$ if it exists, and the union $Q_{Xx}^{n}$\ of the $3$ tiles
with common vertex $x$\ belonging to $\mathcal{P}_{X\Lambda }^{n}$ if they
exist. We are going to prove the following property, which implies the other
statements of Theorem 3.7 and Theorem 3.8:

\noindent For any sets $X,Y$ such that $\mathcal{P}_{X\Lambda },\mathcal{P}%
_{Y\Lambda }$ exist, for any coverings $\mathcal{C},\mathcal{D}$ of $%
\mathcal{P}_{X\Lambda },\mathcal{P}_{Y\Lambda }$ by sets of curves which
satisfy (P), for each $n\in 
%TCIMACRO{\U{2115} }%
%BeginExpansion
\mathbb{N}
%EndExpansion
$ and for each $x\in 
%TCIMACRO{\U{211d} }%
%BeginExpansion
\mathbb{R}
%EndExpansion
^{2}$ such that $Q_{Xx}^{n}$\ exists, each $P_{Yy}^{n+4}$ contains some $%
Q_{Yz}^{n}$\ with $\mathcal{D}\upharpoonright Q_{Yz}^{n}\cong \mathcal{C}%
\upharpoonright Q_{Xx}^{n}$.

By Corollary 3.3, for any $y,z\in 
%TCIMACRO{\U{211d} }%
%BeginExpansion
\mathbb{R}
%EndExpansion
^{2}$ such that $P_{Yy}^{n+2}$ and $P_{Yz}^{n+4}$ exist, each covering of $%
P_{Yz}^{n+4}$ contains copies of the $6$ coverings of $P_{Yy}^{n+2}$.
Consequently, it suffices to show that, for each $x\in 
%TCIMACRO{\U{211d} }%
%BeginExpansion
\mathbb{R}
%EndExpansion
^{2}$ such that $Q_{Xx}^{n}$\ exists, each $P_{Yy}^{n+2}$ contains some $%
Q_{Yz}^{n}$\ such that $\mathcal{C}\upharpoonright Q_{Xx}^{n}$ and $\mathcal{%
D}\upharpoonright Q_{Yz}^{n}$ are equivalent modulo a rotation of angle $%
2k\pi /3$ and/or changing the orientation of all the curves.

We write $\Lambda =(\lambda _{k})_{k\in 
%TCIMACRO{\U{2115} }%
%BeginExpansion
\mathbb{N}
%EndExpansion
^{\ast }}$, $X=(x_{k})_{k\in 
%TCIMACRO{\U{2115} }%
%BeginExpansion
\mathbb{N}
%EndExpansion
}$ and $Y=(y_{k})_{k\in 
%TCIMACRO{\U{2115} }%
%BeginExpansion
\mathbb{N}
%EndExpansion
}$. We have $\Delta _{x_{n}\lambda _{1}\cdots \lambda _{n}}(\mathcal{C}%
)\upharpoonright \Delta _{x_{n}\lambda _{1}\cdots \lambda _{n}}(P)=\Delta
_{x_{n}\lambda _{1}\cdots \lambda _{n}}(\mathcal{C}\upharpoonright P)$ for
each $P\in \mathcal{P}_{X\Lambda }^{n}$ and $\Delta _{y_{n}\lambda
_{1}\cdots \lambda _{n}}(\mathcal{D})\upharpoonright \Delta _{y_{n}\lambda
_{1}\cdots \lambda _{n}}(P)=\Delta _{y_{n}\lambda _{1}\cdots \lambda _{n}}(%
\mathcal{D}\upharpoonright P)$ for each $P\in \mathcal{P}_{Y\Lambda }^{n}$.

For any $x,y\in 
%TCIMACRO{\U{211d} }%
%BeginExpansion
\mathbb{R}
%EndExpansion
^{2}$, $\mathcal{C}\upharpoonright Q_{Xx}^{n}$ and $\mathcal{D}%
\upharpoonright Q_{Yy}^{n}$ are equivalent modulo a translation, or a
rotation of angle $2k\pi /3$, or changing the orientation of the curves, or
a combination of these operations, if and only if the same property is true
for $\Delta _{x_{n}\lambda _{1}\cdots \lambda _{n}}(\mathcal{C}%
)\upharpoonright \Delta _{x_{n}\lambda _{1}\cdots \lambda _{n}}(Q_{Xx}^{n})$
and $\Delta _{y_{n}\lambda _{1}\cdots \lambda _{n}}(\mathcal{D}%
)\upharpoonright \Delta _{y_{n}\lambda _{1}\cdots \lambda _{n}}(Q_{Yy}^{n})$%
. Consequently, it suffices to prove the statement for $n=0$.

Figure 3.5 below shows representatives of the $9$ equivalence classes of
possible configurations for $\mathcal{C}\upharpoonright Q_{Xx}^{0}$, modulo
the operations mentioned just above. We see from Figure 3.3 that, for each
curve $C$ which is a covering of some $P_{Yy}^{2}$, each of these classes is
realized by $C\upharpoonright Q_{Yz}^{0}$ for some $Q_{Yz}^{0}\subset
P_{Yy}^{2}$.~~$\blacksquare $

\begin{center}
\includegraphics[scale=0.60]{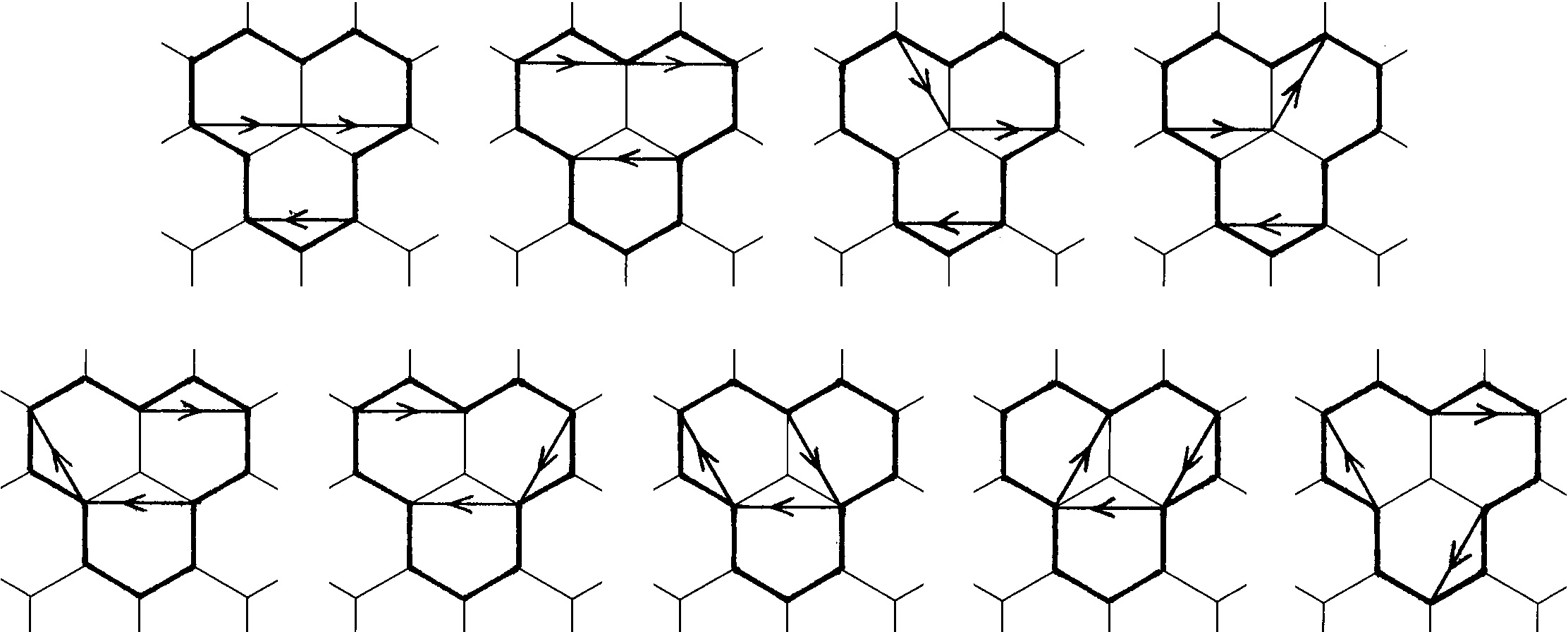}\medskip

Figure 3.5
\end{center}

\bigskip \medskip

\begin{center}
\textbf{References}\medskip
\end{center}

\noindent \lbrack 1] D. Chandler and D.E. Knuth, Number Representations and
Dragon Curves I and II , Journal of Recreational Mathematics, volume 3,
number 2, April 1970, pages 66-81, and number 3, July 1970, pages 133-149.

\noindent \lbrack 2] W. Kuhirun, Simple procedure for evaluating the
impedence matrix of fractal and fractile arrays, Progress in
Electromagnetics Research 14 (2010), 61-70.

\noindent \lbrack 3] B. Mandelbrot, The Fractal Geometry of Nature, W.H.
Freeman and Company, New York, 1983.

\noindent \lbrack 4] F. Oger, Paperfolding sequences, paperfolding curves
and local isomorphism, Hiroshima Math. Journal 42 (2012), 37-75.

\noindent \lbrack 5] F. Oger, The number of paperfolding curves in a
covering of the plane, Hiroshima Math. Journal 47 (2017), 1-14.

\bigskip \medskip

Francis OGER

UFR de Math\'{e}matiques, Universit\'{e} Paris 7

B\^{a}timent Sophie Germain

8 place Aur\'{e}lie Nemours

75013 Paris

France

E-mail: oger@math.univ-paris-diderot.fr

\noindent \vfill

\end{document}